
\def\LaTeX{%
  \let\Begin\begin
  \let\End\end
  \let\salta\relax
  \let\finqui\relax
  \let\futuro\relax}

\def\UK{\def\our{our}\let\sz s}
\def\USA{\def\our{or}\let\sz z}

\UK



\LaTeX

\UK


\salta

\documentclass[twoside,12pt]{article}
\setlength{\textheight}{24cm}
\setlength{\textwidth}{17cm}
\setlength{\oddsidemargin}{2mm}
\setlength{\evensidemargin}{2mm}
\setlength{\topmargin}{-15mm}
\parskip2mm


\usepackage[usenames,dvipsnames]{color}
\usepackage{amsmath}
\usepackage{amsthm}
\usepackage{amssymb}
\usepackage{bbm}
\usepackage[mathcal]{euscript}
\usepackage{graphicx}

\usepackage{hyperref}
\usepackage{enumitem}

\usepackage[utf8]{inputenc}
\usepackage[T1]{fontenc}
\usepackage{mathtools}
\usepackage[thinc]{esdiff}
\newcommand{\dd}{\mathop{}\!\mathrm{d}}
\numberwithin{equation}{section}

\allowdisplaybreaks[4]

%
%


\definecolor{viola}{rgb}{0.3,0,0.7}
\definecolor{ciclamino}{rgb}{0.5,0,0.5}
\definecolor{rosso}{rgb}{0.85,0,0}
\definecolor{verdescuro}{rgb}{0,0.7,0}

\def\an #1{{\color{rosso}#1}}
\def\an #1{#1}
\def\comm #1{{\color{blue}#1}}
\def\comm #1{#1}
\def\modLuca #1{{\color{verdescuro}#1}}
\def\modLuca #1{#1}





\bibliographystyle{plain}


%
\newtheorem{theorem}{Theorem}[section]
\newtheorem{remark}[theorem]{Remark}

\newtheorem{definition}[theorem]{Definition}

\finqui

\def\Bcenter{\Begin{center}}
\def\Ecenter{\End{center}}
\let\non\nonumber




\def\step #1 \par{\medskip\noindent{\bf #1.}\quad}


\def\rhs{right-hand side}



\def\multibold #1{\def\arg{#1}%
  \ifx\arg\pto \let\next\relax
  \else
  \def\next{\expandafter
    \def\csname #1#1#1\endcsname{{\boldsymbol #1}}%
    \multibold}%
  \fi \next}

\def\pto{.}

\def\multical #1{\def\arg{#1}%
  \ifx\arg\pto \let\next\relax
  \else
  \def\next{\expandafter
    \def\csname #1#1\endcsname{{\cal #1}}%
    \multical}%
  \fi \next}


\def\multimathop #1 {\def\arg{#1}%
  \ifx\arg\pto \let\next\relax
  \else
  \def\next{\expandafter
    \def\csname #1\endcsname{\mathop{\rm #1}\nolimits}%
    \multimathop}%
  \fi \next}

\multibold
qwertyuiopasdfghjklzxcvbnmQWERTYUIOPASDFGHJKLZXCVBNM.

\multical
QWERTYUIOPASDFGHJKLZXCVBNM.


\multimathop
diag dist div dom mean meas sign supp .


\def\Accorpa #1#2 #3 {\gdef #1{\eqref{#2}-\eqref{#3}}%
  \wlog{}\wlog{\string #1 -> #2 - #3}\wlog{}}
\def\Accorparef #1#2 #3 {\gdef #1{\ref{#2}-\ref{#3}}%
  \wlog{}\wlog{\string #1 -> #2 - #3}\wlog{}}


\def\<#1>{\mathopen\langle #1\mathclose\rangle}

\def\[#1]{\mathopen\langle\!\langle #1\mathclose\rangle\!\rangle}

\def\checkmmode #1{\relax\ifmmode\hbox{#1}\else{#1}\fi}


\def\erre{{\mathbb{R}}}




\def\genspazio #1#2#3#4#5{#1^{#2}(#5,#4;#3)}
\def\spazio #1#2#3{\genspazio {#1}{#2}{#3}T0}

\def\L {\spazio L}

\def\C #1#2{{\cal C}^{#1}([0,T];#2)}


\def\Lx #1{L^{#1}(\Omega)}



\let\theta\vartheta

\let\eps\varepsilon
\let\ph\varphi

\let\TeXchi\chi                         
\newbox\chibox
\setbox0 \hbox{\mathsurround0pt $\TeXchi$}
\setbox\chibox \hbox{\raise\dp0 \box 0 }
\def\chi{\copy\chibox}




\def\hh{{\mathbbm{h}}}
\def\ov #1{{\overline{#1}}}

\def\0{{\boldsymbol{0}}}

\def\hhu{\hh_1}
\def\hhd{\hh_2}

\Begin{document}


%
\title{On a non-local phase-field model for tum\our\ growth with single-well Lennard--Jones potential}

\date{}
\author{}
\maketitle
\Bcenter
\vskip-1cm
{\large\sc Maurizio Grasselli$^{(1)}$}\\
{\normalsize e-mail: {\tt maurizio.grasselli@polimi.it}}\\[.25cm]
{\large\sc Luca Melzi$^{(2)}$}\\
{\normalsize e-mail: {\tt l.melzi24@imperial.ac.uk}}\\[.25cm]
{\large\sc Andrea Signori$^{(1),(3)}$}\\
{\normalsize e-mail: {\tt andrea.signori@polimi.it}}\\[.5cm]
$^{(1)}$
{\small Dipartimento di Matematica, Politecnico di Milano}\\
{\small via E. Bonardi 9, 20133 Milano, Italy}
\\[.3cm]
$^{(2)}$
{\small Department of Mathematics, Imperial College London}\\
{\small London, SW7 2AZ, UK}
\\[.3cm]
$^{({3})}$
{\small Alexander von Humboldt Research Fellow}
\\[1cm]\date{}

\Ecenter

\Begin{abstract}
\noindent
In the present work, we develop a comprehensive and rigorous analytical framework for a non-local phase-field model that describes tum\our\ growth dynamics. The model is derived by coupling a non-local Cahn--Hilliard equation with a parabolic reaction-diffusion equation, which accounts for both phase segregation and nutrient diffusion.
Previous studies have only considered symmetric potentials for similar models. However, in the biological context of cell-to-cell adhesion, single-well potentials, like the so-called Lennard--Jones potential, seem physically more appropriate. The Cahn--Hilliard equation with this kind of potential has already been analysed. Here, we take a step forward and consider a more refined model.
First, we analyse the model with a viscous relaxation term in the chemical potential and subject to suitable initial and boundary conditions. We prove the existence of solutions, a separation property for the phase variable, and a continuous dependence estimate with respect to the initial data.
Finally, via an asymptotic analysis, we recover the existence of a weak solution to the initial and boundary value problem without viscosity, provided that the chemotactic sensitivity is small enough.

\vskip3mm

\noindent {\bf Keywords:} Tum\our\ growth, phase-field model, Cahn--Hilliard equation, non-local interactions, Lennard--Jones potential, viscous relaxation, weak and strong solutions, continuous dependence on the initial datum.


\End{abstract}

\salta
\pagestyle{myheadings}
\newcommand\testopari{\sc Grasselli \ -- \ Melzi \ -- \ Signori}
\newcommand\testodispari{{\sc tum\our\ growth model with single-well potential}}
\markboth{\testodispari}{\testopari}
\finqui

\section{Introduction}
Cancer is a leading cause of premature death worldwide, and mathematical modelling plays a crucial role in understanding and combating it. In recent decades, phase-field models have been proposed to describe tum\our\ growth dynamics, where the tum\our\ is modelled using an order parameter representing the local concentration of tum\our\ cells, see, e.g., \cite{Bellomo2008,Cristini2010,Friedman2007}. 
These models have been extensively studied, with many examples in the literature, such as \cite{Garcke2016}, and the references therein.
These are capable of accounting for complex behavi\our s\ like topological changes in tum\our\ regions and interactions with surrounding tissue, such as nutrient consumption. 
Although the phase segregation described by the local Cahn--Hilliard equation is widely accepted, it fails to capture long-range cell-to-cell and cell-to-matrix adhesion, crucial for tum\our\ growth (see, for instance, \cite{Wise2008}). Ignoring long-range interactions is a significant drawback, particularly in modelling tum\our -cell invasion and metastasis formation, processes driven by long-range competition.
To address this, non-local Cahn--Hilliard equations have been used, incorporating interactions at a distance through modified energy functionals (see \cite{Giacomin1996,Giacomin1997,Giacomin1998}).
These models, which have been rigorously validated, are crucial for accurately capturing processes such as cell invasion, which involve long-range competition. Notably, while cell adhesion is inherently a non-local phenomenon in space, the chemotaxis mechanism, in contrast, is of a local nature, as emphasi\sz ed in \cite{sherratt2009boundedness}.
Here we consider a variant of the non-local Cahn--Hilliard equations that is characterised by a single-well Lennard--Jones potential density. This choice better represents cell-to-cell adhesion and aims to contribute to the growing body of models that can aid in understanding and treating cancer (see \cite{Agosti2017,Agosti2018bis,Agosti2024} and references therein).

In this contribution, we consider a class of models inspired by \cite{Garcke2016}, which has been extended to the non-local setting in \cite{Scarpa2021}.
Although the authors of \cite{Scarpa2021} focus on regular or singular double-well symmetric potentials, this study extends their analysis to the singular single-well potential, which proves to be relevant from a biological point of view (cf. \cite{Agosti2017}).

Let $\Omega\subset\mathbb{R}^d$, $d \in \{2,3\}$, be a bounded smooth domain and set a final time $T>0$. The problem we want to analyse is the following
\begin{alignat}{2}
    &\partial_t\varphi-\Delta\mu=S(\varphi,\sigma)\qquad&&\text{in }Q := \Omega \times (0,T), \label{eq:varphi} \\
    &\mu=\tau\partial_t\varphi+ \eps a\varphi-\eps J*\varphi+ \eps^{-1}F'(\varphi)-\chi\sigma\qquad&&\text{in }Q, \label{eq:mu} \\
    &\partial_t\sigma-\Delta\sigma+B(\sigma-\sigma_S)+C\sigma {\hhd}(\varphi)=0\qquad&&\text{in }Q, \label{eq:sigma} \\
    &\partial_{\boldsymbol{n}}\mu=\partial_{\boldsymbol{n}}\sigma=0\qquad&&\text{on }\Sigma:= \partial\Omega\times(0,T), \label{eq:bcs} \\
    &\varphi(0)=\varphi_0, \ \sigma(0)=\sigma_0 \qquad&&\text{in }\Omega.
    \label{eq:ics}
\end{alignat}
\Accorpa\Sys {eq:varphi} {eq:ics}
Here, $\varphi$ represents the phase variable, which denotes the tum\our\ concentration in terms of volume, normali\sz ed between zero and one. The parameter $\eps > 0$ defines the thickness of the interface between pure phases. With this choice, the level sets $\{\varphi = 0\}$ and $\{\varphi = 1\}$ correspond to regions of healthy and tum\our\ cells, respectively. For sufficiently large times, we expect that most of the domain is occupied by these two pure phases, separated by the thin diffuse interface $\{0 < \varphi < 1\}$. Since it does not affect our analysis, we set $\eps = 1$ for convenience from now on. 
Next, $\mu$ represents the chemical potential, and $\sigma$ is the concentration of a nutrient species that satisfies the reaction-diffusion equation \eqref{eq:sigma}. Since both $\varphi$ and $\sigma$ represent concentrations, it is natural to assume that they are within the interval $[0, 1]$. This will be rigorously established later using the singular nature of the potential in the model. Here $\mu$ is given \modLuca{by} the variational derivative of a suitable Helmholtz free energy, but it also contains
an additional viscosity term $\tau\partial_t\varphi$, $\tau \in (0,1]$
(see \cite{Novick-Cohen1988}). We recall that the non-local free energy is given by
\begin{align*}
	E(\varphi) &= \frac{1}{4}\int_{\Omega \times \Omega} J(x - y) |\varphi(x) - \varphi(y)|^2 \dd x \dd y + \int_\Omega F(\varphi),
\end{align*}
where $F$ is a suitable single-well or double-well potential, and $J$ is a spatially symmetric convolution kernel suitably peaked around zero. Specifically, we define:
\begin{align*}
	(J * v)(x) &= \int_\Omega J(x - y) v(y) \dd y\quad\forall x \in \Omega, \forall v: \Omega \to \mathbb{R} \,\text{ measurable},
\end{align*}
and $a(x) := (J * 1)(x)$ for every $x \in \Omega$.

A common choice for the potential $F$ is the double-well form $F(r) = \frac 14 r^2(r-1)^2$, which appears in many models in the literature, usually normali\sz ed between $-1$ and $1$. However, this choice does not ensure that the phase field \modLuca{variable} takes its values in the physical range (say, $[0,1]$). Indeed, a more reasonable choice is the logarithmic potential (see, e.g., \cite{Scarpa2021}). In contrast, this study focuses on a single-well potential, specifically the so-called Lennard\an{--}Jones potential.
This potential captures the asymmetry in cell-cell interactions, which are attractive at moderate tum\our\ concentrations when $\varphi < \bar{\varphi}$, for some $\bar{\varphi} \in (0,1)$, and repulsive at higher concentrations when $\varphi > \bar{\varphi}$, respectively. The Lennard\an{--}Jones potential is defined as
\begin{align*}
    F(r)=\begin{cases}
        -(1-\bar{\varphi})\log(1-r)-\frac{r^3}{3}-(1-\bar\varphi)(\frac{r^2}{2}+r)&\text{ if }r\in[0,1),\\
        +\infty&\text{ otherwise},
    \end{cases}
\end{align*}
where $F'(\bar\varphi)=0$.

Furthermore, in equation \eqref{eq:sigma}, $B$ and $C$ are non-negative constants, ${\hhd}$ is a real-valued function, and $\sigma_S$ represents a pre-existing nutrient concentration.
Finally, $S(\ph,\sigma)$ in \eqref{eq:varphi} acts as a mass source, as it is directly linked to the mass conservation equation and follows a generali\sz ed Cahn--Hilliard--Oono form:
\begin{align*} 
	S(\varphi, \sigma) &= -m \varphi + {\hhu}(\varphi, \sigma), 
\end{align*}
where $m > 0$ and ${\hhu}: \mathbb{R}^2 \to \mathbb{R}$ a suitable uniformly bounded function.
 When ${\hhu} \equiv c \in (0,m)$, we recover the so-called Cahn--Hilliard--Oono equation for $\varphi$, see, e.g., \cite{fukao2024cahn, GGM, miranville2016cahn}, and the references therein.
Finally, the constant $\chi > 0$ in \eqref{eq:mu} represents the chemotactic sensitivity, modelling the tendency of tum\our\ cells to migrate toward areas of higher nutrient concentration. For further details on the mathematical modelling of chemotaxis and related mechanisms, such as active transport, we refer to \cite{Garcke2017,Garcke2018, Garcke2016, he2022viscous, Rocca2023}. Although active transport is excluded from our study, it is considered in some of the results in \cite{Scarpa2021} which, however, require a further relaxation of the chemical potential.

It is also worth mentioning that the non-local Cahn--Hilliard equation alone,
i.e.\an{,} \eqref{eq:varphi}-\eqref{eq:mu} without source has been analysed in \cite{Agosti2024}, but taking a suitable degenerate mobility. This choice seems
physically more appropriate (see also \cite{Agosti2017} and references therein) and it will be the subject of further investigations.

Without aiming to be exhaustive, we recall here below a selection of works from the literature that are related to the model under investigation.  
In addition to the studies mentioned above, for local models involving double-well potentials, we refer to \cite{Colli2015,  
 Garcke2017, Garcke2018,Garcke2016}. Non-local models have been explored in \cite{frigeri2022strong,Scarpa2021}, while single-well models have been considered in \cite{Agosti2017,Agosti2018bis, Agosti2024, PP}.  
Several works have also addressed models that incorporate velocity effects, including \cite{Colli2023, Ebenbeck2019, Ebenbeck2021, Garcke2016bis, Garcke2018bis, Jiang2015}. Multiphase models have been widely studied, with notable contributions such as \cite{agosti2024analysis, frigeri2017multi, knopf2022existence}.  
Other relevant mechanisms, such as elasticity, have been investigated in \cite{garcke2022viscoelastic, garcke2021phase}.  
Finally, we highlight applications to optimal control problems (see \cite{Colli2021, Colli2022bis, Gilardi2023bis, Signori2020bis, signoricontrol}, cf. also  \cite{fornoni2024maximal, fornoni2024optimal, Rocca2021} for non-local Cahn--Hilliard equations).


The paper is organi\sz ed as follows. In Section \ref{sec:Setting and main results}, we introduce the notation and key properties, list our basic assumptions, and state our main results.  
In Section \ref{sec:analysis_with_tau>0}, we consider the viscous case $\tau > 0$, and we prove the existence of weak solutions. Further results are proven in Section \ref{sec:moreanalysis_with_tau>0}. First, we demonstrate the existence of more regular solutions and the fact that $\varphi$ is strictly separated from the pure phase $1$ if the initial datum does (cf. \cite{Agosti2024}). Then, we establish a continuous dependence estimate for regular and strictly separated solutions, ensuring their uniqueness.
Finally, in Section \ref{sec:asymptotics_as_tau_to_0}, we analy\sz e the asymptotic behavi\our\ as $\tau \to 0^+$, showing that problem \Sys\ with $\tau = 0$ itself admits a weak solution.

\section{Setting and main results} \label{sec:Setting and main results}
Before presenting the results, let us introduce some notation that will be used throughout this work.
We consider a bounded smooth domain $\Omega \subset \mathbb{R}^d$, with $d \in \{2,3\}$, and denote by $|\Omega|$ its Lebesgue measure. We also define $\Gamma := \partial \Omega$ as the boundary of $\Omega$.
Fixing $T>0$, we set
\begin{align*}
    Q:=\Omega\times(0,T),\quad\Sigma:=\Gamma \times (0,T),\quad Q_t:=\Omega\times(0,t)\quad\forall t\in(0,T).
\end{align*}
Then, we set $H:=L^2(\Omega)$ and $V:=H^1(\Omega)$, so that $(V,H,V')$ is an Hilbert triplet, with $V'$ denoting the topological dual of $V$ and $\<\cdot,\cdot>$ the associated duality pairing between $V'$ and $V$. The scalar products are denoted by $(\cdot,\cdot)_V$, $(\cdot,\cdot)$ and $(\cdot,\cdot)_{V'}$ respectively, and the norms $\|\cdot\|_V$, $\|\cdot\|$ and $\|\cdot\|_{V'}$ are defined accordingly.
As usual, identifying $H$ with its dual $H'$, we can write $V\hookrightarrow H\simeq H'\hookrightarrow V'$, where all the inclusions are dense, continuous and compact.
Furthermore, we consider
\begin{align*}
    W:=\{v\in H^2(\Omega)\mid\partial_{\boldsymbol{n}}v=0\quad\text{a.e. on }\Gamma\},
\end{align*}
with its natural norm denoted by $\|\cdot\|_W$.
We have that $W\subset V\subset H\subset V'\subset W'$, where all the inclusions are dense, continuous and compact.
We define the generali\sz ed mean value of  $v\in V'$ as $v_\Omega:=\frac{1}{|\Omega|}\langle v,1\rangle$.
Moreover, recall that the variational operator $-\Delta:V\to V'$, given by
\begin{align*}
    \langle-\Delta u,v\rangle:=\int_\Omega\nabla u\cdot\nabla v\quad\forall u,v\in V,
\end{align*}
is such \modLuca{a way} that its inverse $\mathcal{N}:=(-\Delta)^{-1}$ is an isomorphism from $V'_0:=\{v\in V'\mid v_\Omega=0\}$ to $V_0:=\{v\in V\mid v_\Omega=0\}$, and satisfies
\begin{alignat}{2}
    \langle-\Delta u,\mathcal{N}v\rangle&=\langle v,u\rangle\quad&&\forall u\in V_0,\forall v\in V_0', \label{eq:property_of_N_1} \\
    \langle u,\mathcal{N}v\rangle&=\langle v,\mathcal{N}u\rangle=(u,v)_{V'}\quad&&\forall u,v\in V'_0, \label{eq:property_of_N_2} \\
    \langle\partial_t v(t),\mathcal{N}v(t)\rangle&=\frac{1}{2}\diff{}{t}\|v(t)\|^2_{V'}\quad&&\text{for a.e. }t\in(0,T),\forall v\in H^1(0,T;V_0'). \label{eq:property_of_N_3}
\end{alignat}
Finally, throughout the paper we will use the letter $M$ to denote a generic positive constant that depends only on the data of the problem \eqref{eq:varphi}-\eqref{eq:ics}, and its value may change from line to line.
Possibly, we will write $M_\eta$ to specify the dependence on some positive parameter $\eta$.

\subsection{Statements of the results}

To begin with, we make the following structural assumptions about the data.
\begin{enumerate}[label={\bf (A\arabic{*})}, ref={\bf (A\arabic{*})}]

	\item \label{eq:assumption1} 		
	${\hhu}\in \mathcal{C}^0(\mathbb{R}^2)\cap L^\infty(\mathbb{R}^2)\,\,\text{with }K:=\|{\hhu}\|_{L^\infty(\mathbb{R}^2)}<m, $
	and 
	${\hhu}(r,s)\geq0\quad\forall(r,s)\in\mathbb{R}^2.$ 
	
\item \label{eq:assumption2}
	${\hhd}\in \mathcal{C}^0(\mathbb{R})\cap L^\infty(\mathbb{R})$.
	
	\item  \label{eq:assumption3}  $\sigma_S\in L^\infty(Q)$ with $\sigma_S(x,t)\in[0,1]\quad\text{for a.e. }(x,t)\in Q.$ 
	
		\item \label{ass:LJpot} $F$ is the Lennard--Jones potential density given by 	\begin{align}
    \label{def:LJ}
    F(r)=\begin{cases}
        -h\log(1-r)-\frac{r^3}{3}-h(\frac{r^2}{2}+r)&\text{ if }r\in[0,1),\\
        +\infty&\text{ otherwise},
    \end{cases}
\end{align}
with $\bar{\varphi} \in (0,1)$ and $h:=1-\bar{\varphi}$.
Moreover, we split the potential into a convex, singular component $F_1$ and a non-convex, regular component $F_2$, in such a way that $F=F_1+F_2$, as follows:
\begin{align*}
    F_1(r)&:=\begin{cases}
        -h\log(1-r)&\text{ if }r\in[0,1),\\
        +\infty&\text{ otherwise},
    \end{cases} \qquad F_2(r):=\begin{cases}
        -\frac{r^3}{3}-h(\frac{r^2}{2}+r)&\text{ if }r<1,\\
        +\infty&\text{ otherwise}.
    \end{cases}
\end{align*}

    
    \item \label{ass:kernel} The spatial convolution kernel $J \in W^{1,1}_{\rm loc}(\mathbb{R}^d)$ and satisfies
\begin{align*}
    J(x) = J(-x), \quad \text{for almost all }x \in \erre^d.
\end{align*}
Furthermore, we postulate that
\begin{align*}
    a_* &:= \inf_{x \in \Omega} \int_\Omega J(x-y) \, dy = \inf_{x \in \Omega} a(x)  \geq2+h-3h^\frac{1}{3}=:c_0>0, 
    \\ 
    a^* & := \sup_{x \in \Omega} \int_\Omega |J(x-y)| \, dy < +\infty, \quad
    b^* := \sup_{x \in \Omega} \int_\Omega |\nabla J(x-y)| \, dy < +\infty.
\end{align*}
\end{enumerate}
We observe that the strict positivity of $a^*$ entails
\begin{equation}
 F''(r)+a_*\geq c_0\quad\forall r\in\mathbb{R}.\label{eq:property_of_F_second_0}
\end{equation}
This is necessary in order to ensure the well-posedness of the non-local Cahn--Hilliard equation.

Before presenting the list of our results, we begin by providing the following definition.
\begin{definition}[Weak solution]
    \label{def:weak_sol}
    Let $\tau$ be a nonnegative constant. The triplet $(\varphi,\mu,\sigma)$ is called a weak solution to problem \eqref{eq:varphi}-\eqref{eq:ics} whenever
    \begin{align*}
      &   \varphi\in H^1(0,T;V')\cap L^2(0,T;V),\quad
        \mu\in L^2(0,T;V),
        \\
        & \sigma\in H^1(0,T;V')\cap L^2(0,T;V),
    \end{align*}
    and
    \begin{align*}
      &  \langle \partial_t\varphi,\xi \rangle + \int_\Omega\nabla\mu\cdot\nabla\xi=\int_\Omega S(\varphi,\sigma)\xi,\\
      &  \langle \partial_t\sigma,\xi \rangle+\int_\Omega\nabla\sigma\cdot\nabla\xi+\int_\Omega \biggl(B(\sigma-\sigma_S)+C\sigma {\hhd}(\varphi)\biggr)\xi=0,
    \end{align*}
    for all $\xi \in V$, almost everywhere in $(0,T)$, and
    \begin{alignat*}{2}
       & \mu=\tau\partial_t\varphi+a\varphi-J\ast\varphi+F'(\varphi)-\chi\sigma &&\quad \text{a.e. in }Q,\\
       & \varphi(0)=\varphi_0,\quad\sigma(0)=\sigma_0 &&\quad \text{a.e. in }\Omega.
    \end{alignat*}
\end{definition}
Our first result deals with the existence of weak solutions to problem \eqref{eq:varphi}-\eqref{eq:ics} for $\tau>0$.

\begin{theorem}
[Existence of a weak solution]
\label{thm:existence_weak_tau>0}
    Let $\tau\in (0,1]$ be fixed. Assume that \ref{eq:assumption1}-\ref{ass:kernel} hold.
    Moreover, let the initial data fulfil
    \begin{align}
        \varphi_0\in V,\quad(\varphi_0)_\Omega\in(0,1),\quad F(\varphi_0)\in L^1(\Omega),\quad\sigma_0\in H.
        \label{eq:assumptions_ics}
    \end{align}
    Then, there exists a triplet $(\varphi,\mu,\sigma)$ that is a weak solution to \eqref{eq:varphi}-\eqref{eq:ics} in the sense of Definition \ref{def:weak_sol}.
    Moreover, it holds that
    \begin{align}
        & \varphi\in H^1(0,T;H)\cap L^\infty(0,T;V)\cap L^\infty(Q), \label{eq:weak_regularity_varphi} \\
        & \varphi\in[0,1)\quad\text{a.e. in }Q. \label{eq:MP_varphi} 
    \end{align}
     If we also suppose
    \begin{align}
        &\sigma_0(x)\in[0,1]\quad\text{for a.e. }x\in\Omega, \label{eq:sigma0_in01} \\
        &{\hhd}(r)\geq0\quad\forall r\in\mathbb{R}, \label{eq:h2_nonnegative}
    \end{align}
    then, $\sigma\in L^\infty(Q)$ with
    \begin{equation}
        \sigma\in[0,1]\quad\text{a.e. in }Q. \label{eq:MP_sigma}
    \end{equation}   
\end{theorem}

\begin{remark}
It is worth noticing that \eqref{eq:MP_varphi} is related to the fact that $F$ is singular at $1$ (see \eqref{def:LJ}).  
\end{remark}

Our second result deals with more regular (i.e., strong) solutions.
\begin{theorem} 
[Existence of a strong solution]
\label{thm:strong_solutions_tau>0}
Let $\tau\in (0,1]$ be fixed. Assume that \ref{eq:assumption1}-\ref{ass:kernel} and  \eqref{eq:assumptions_ics} hold.
Suppose further that
    \begin{align}
        & F'(\varphi_0)\in V, \quad \sigma_0\in V, \label{eq:assumptions_ics_enhanced} \\
        &{\hhu}\text{ is Lipschitz continuous in }\mathbb{R}^2. \label{eq:grad_h1_bdd}
    \end{align}
    Then, a triplet $(\varphi,\mu,\sigma)$ given by Theorem \ref{thm:existence_weak_tau>0} also satisfies 
    \begin{align}
        \varphi&\in W^{1,\infty}(0,T;H) \cap \mathcal{C}^0
        ([0,T];V), \label{eq:additional_reg_varphi} \\
        \mu&\in L^\infty(0,T;V), \label{eq:additional_reg_mu} \\
        \sigma&\in H^1(0,T;H)\cap L^2(0,T;W). \label{eq:regularity_sigma_bis}
    \end{align}

\end{theorem}

\begin{remark}
Observe that a solution given by Theorem \ref{thm:strong_solutions_tau>0}
solves problem \eqref{eq:varphi}-\eqref{eq:ics} almost everywhere. Thus, it is a strong solution.
Moreover, due to \eqref{eq:regularity_sigma_bis}, $\sigma \in \mathcal{C}^0([0,T];V)$.
It is worth mentioning that the existence of a weak (or strong) solution actually requires only $\sigma_S\in L^2(Q)$ (cf. Assumption \ref{eq:assumption3}).
However, the meaning of $\sigma$ requires its boundedness.
\end{remark}

\begin{remark}
Let $(\varphi,\mu,\sigma)$ be a weak solution and suppose that \eqref{eq:grad_h1_bdd} holds. Thanks to \eqref{eq:weak_regularity_varphi} and the fact that $\sigma\in L^2(0,T;V)$, for any fixed $t\in (0,T)$, we can find $t^* \in (0,t)$ such that 
$$
F'(\varphi(t^*))\in V, \quad \sigma(t^*)\in V.
$$
Thus, recalling \eqref{eq:assumptions_ics_enhanced}, we see that $(\varphi,\mu,\sigma)$
is a strong solution in $(t^*,T)$.
\end{remark}

Next, we exploit this additional regularity to prove that $\varphi$ stays uniformly away from $1$, that is, $\varphi$ satisfies the strict separation property, provided that $\varphi_0$ does.

\begin{theorem}[Strict separation property] \label{thm:separation_property_varphi_tau>0}
Let $\tau \in (0,1]$ be fixed. 
Assume that \ref{eq:assumption1}-\ref{ass:kernel}, \eqref{eq:assumptions_ics}, 
\eqref{eq:sigma0_in01}-\eqref{eq:h2_nonnegative}, \eqref{eq:assumptions_ics_enhanced}-\eqref{eq:grad_h1_bdd}  
hold. If 
$(\varphi,\mu,\sigma)$
is a strong solution 
given by Theorem \ref{thm:strong_solutions_tau>0}
and there exists $\bar{\delta}\in(0,1)$ such that
    \begin{align}
        0\leq\varphi_0\leq1-\bar{\delta}\quad\text{a.e. in }\Omega,
        \label{eq:separation_for_ic}
    \end{align}
    then there exists $\delta\in(0,\bar{\delta}]$, possibly dependent on $\tau$, such that
    \begin{align}
        0\leq\varphi\leq1-\delta\quad\text{a.e. in }Q.
\label{eq:separation_property_varphi}
    \end{align}
\end{theorem}

Next, we state the continuous dependence of the solution on the initial conditions.
\begin{theorem}[Continuous dependence on the initial data] \label{thm:continuous_dependence_tau>0}
    Let $\tau \in (0,1]$ be fixed.
    Assume that \ref{eq:assumption1}-\ref{ass:kernel} hold.
    Let $(\varphi_0^j,\sigma_0^j)$, $j=1,2$, satisfy the assumptions of Theorem \ref{thm:separation_property_varphi_tau>0} which are supposed to hold also for $\hhu$ and $\hhd$. Suppose further that ${\hhd}$ is Lipschitz continuous on $\mathbb{R}$. Let $(\varphi_j,\mu_j,\sigma_j)$ be a strong solution to \eqref{eq:varphi}-\eqref{eq:ics} associated with $(\varphi_0^j,\sigma_0^j)$, $j=1,2$. Then, the following estimate holds
    \begin{align}
        \non 
        & \|(\varphi_1-(\varphi_1)_\Omega)-(\varphi_2-(\varphi_2)_\Omega)\|_{L^\infty(0,T;V')}
        +\|\varphi_1-\varphi_2\|_{\mathcal{C}^0([0,T];H)}+\|\sigma_1-\sigma_2\|_{\mathcal{C}^0([0,T];H)\cap L^2(0,T;V)}\\
        & \quad \leq M_\tau\bigl(\|(\varphi_0^1-(\varphi_0^1)_\Omega)-(\varphi_0^2-(\varphi_0^2)_\Omega)\|_{V'}+\|\varphi_0^1-\varphi_0^2\|+\|\sigma_0^1-\sigma_0^2\|\bigr),
        	\label{cd:est}
    \end{align}
    for some constant $M_\tau>0$ that possibly depends on $\tau$.
\end{theorem}
As a consequence, we infer the uniqueness of strictly separated  (strong) solutions.

Finally, we state a result on the asymptotic analysis as $\tau\to 0^+$, which gives the existence of a weak solution to problem \Sys\ when $\tau=0$.
\begin{theorem}[Asymptotics as $\tau\to0^+$] \label{thm:asymptotic_tau_to_0}
    Let \ref{eq:assumption1}-\ref{ass:kernel} hold. Consider two sequences of initial conditions $\{\varphi_{0,\tau}\}_\tau,\{\sigma_{0,\tau}\}_\tau\subset H$ satisfying
    \begin{align}
       & \varphi_{0,\tau}\overset{H}{\to}\varphi_{0},\quad\sigma_{0,\tau}\overset{H}{\to}\sigma_{0},
        \label{eq:lastthm_ics1}
    \end{align}
    as $\tau\to 0^+$, for some $\varphi_0,\sigma_0\in H$ such that $F(\varphi_0)\in L^1(\Omega)$.
    Suppose also that
    \begin{align}
       & (\varphi_{0,\tau})_\Omega\in(0,1)\quad\forall\tau\in(0,1),\\
      &  \chi<\min\bigl\{\sqrt{\tfrac{c_0}{2}},1\bigr\}, \label{eq:tau_critic_0} \\
   &     \|F(\varphi_{0,\tau})\|_{L^1(\Omega)}+\tau^\frac{1}{2}\|\varphi_{0,\tau}\|_V\leq M_0\quad\forall\tau\in(0,1), \label{eq:lastthm_ics2}
    \end{align}
    for some positive constant $M_0$ that is independent of $\tau$.
    For each $\tau\in(0,1)$, denote by $(\varphi_\tau,\mu_\tau,\sigma_\tau)$ a triplet given by Theorem \ref{thm:existence_weak_tau>0}, with $(\varphi_{0,\tau},\sigma_{0,\tau})$ as initial conditions.
    Then, there exists a triplet $(\varphi,\mu,\sigma)$ that is a weak solution to \eqref{eq:varphi}-\eqref{eq:ics} with $\tau=0$ in the sense of Definition \ref{def:weak_sol}, and satisfies
    \begin{align*}
        \varphi&\in H^1(0,T;V')\cap\mathcal{C}^0([0,T];H)\cap L^2(0,T;V), \\
        \mu&\in L^2(0,T;V), \\
        \sigma&\in H^1(0,T;V')\cap\mathcal{C}^0([0,T];H)\cap L^2(0,T;V), \\
        \varphi&\in L^\infty(Q),\quad\text{with }\varphi\in[0,1)\quad\text{a.e. in }Q.
    \end{align*}
    Moreover, as $\tau\to0^+$, up to subsequences,
    \begin{alignat*}{2}
       \varphi_{\tau}& \overset{*}{\rightharpoonup}\varphi\quad\text{in }L^\infty(Q)\cap H^1(0,T;V')\cap L^2(0,T;V),\quad &&\varphi_{\tau}\to\varphi\quad\text{in }L^2(0,T;H), \\
       \mu_{\tau}& \rightharpoonup\mu\quad\text{in }L^2(0,T;V),\quad &&\mu_{\tau}\to\mu\quad\text{in }L^2(0,T;H), \\
        \sigma_{\tau}& \overset{*}{\rightharpoonup}\sigma\quad\text{in }H^1(0,T;V')\cap L^\infty(0,T;H)\cap L^2(0,T;V),\quad&&\sigma_{\tau}\to\sigma\quad\text{in }\mathcal{C}^0([0,T];V')\cap L^2(0,T;H), \\
        \tau\varphi_{\tau}& \to0\quad\text{in }H^1(0,T;H)\cap L^\infty(0,T;V). &&
    \end{alignat*}
\end{theorem}

\section{ Proof of Theorem \ref{thm:existence_weak_tau>0}} \label{sec:analysis_with_tau>0}
This section is devoted to the analysis of problem \Sys\ in the viscous case (i.e., $\tau\in (0,1]$). After a  preliminary subsection on the evolution of the spatial average of $\varphi$, the remaining three subsections are devoted to the proof of existence of a weak solution through a Faedo--Galerkin scheme combined with a suitable approximation of the nonlinear potential $F$. 

\subsection{Controlling the spatial average} \label{subsec:a_priori_estimates}
Let $(\varphi,\mu,\sigma)$ be a solution to the system \eqref{eq:varphi}-\eqref{eq:ics}. Let us test \eqref{eq:varphi} by $|\Omega|^{-1}$.
In view of \eqref{eq:bcs}, thanks to the Divergence Theorem, we infer that
\begin{align*}
    \diff{}{t}\varphi_\Omega+m\varphi_\Omega=\frac{1}{|\Omega|}\int_\Omega {\hhu}(\varphi,\sigma).
\end{align*}
Assuming \ref{eq:assumption1}  and  defining $y:=\varphi_\Omega$, we end up with 
\begin{align*}
    0\leq y'+my\leq K,
\end{align*}
so that
\begin{align*}
    y(0)e^{-mt}\leq y(t)\leq y(0)e^{-mt}+(1-e^{-mt})\frac{K}{m}.
\end{align*}
Recalling \eqref{eq:assumptions_ics}, we have $y(0)=(\varphi_0)_\Omega\in(0,1)$.
Using \ref{eq:assumption1} once more, we realise that there exists $\delta\in(0,\frac{1}{2})$ such that
\begin{align}
    \delta\leq\varphi_\Omega(t)\leq1-\delta\quad\forall t\in[0,T].
    \label{eq:separation_for_mean}
\end{align}
For instance, we can select
\begin{align*}
    \delta:=\min\bigl\{\tfrac{1}{4},(\varphi_0)_\Omega e^{-mT},1-\max\bigl\{\tfrac{K}{m},(\varphi_0)_\Omega\bigr\}\bigr\}.
\end{align*}

\subsection{The approximating problem}
\label{subsec:The finite dimensional regularized problem}
Here we first approximate $F$ (see \eqref{def:LJ}) and then
use a Galerkin discretisation.

Regarding the former, for any $\lambda\in (0,1)$, we set
\begin{align}
    F_\lambda(r)&=F_{1,\lambda}(r)+\bar{F_2}(r), \label{eq:regularized_potential} \\
    F_{1,\lambda}(r)&=
    \begin{cases}
    -\frac{r^3}{\lambda}+\frac{h}{2}r^2+hr&\text{for }r<0,\\
    F_1(r)&\text{for }0\leq r< 1-\lambda, \\
    -h\log\lambda+\frac{3}{2}h-\frac{2}{\lambda}h(1-r)+\frac{h}{2\lambda^2}(1-r)^2&\text{for} \ r\geq\ 1-\lambda,
    \end{cases} \label{eq:regularized_potential1} \\
    \bar{F_2}(r)&=
    \begin{cases}
    F_2(r)&\text{for }r < 1, \\
    -\frac{1}{3}-\frac{3}{2}h+(-1-2h)(r-1)+\frac{1}{2}(-2-h)(r-1)^2&\text{for} \ r\geq 1. \label{eq:regularized_potential2}
\end{cases}
\end{align}
In Figure \ref{fig:01}, we report a plot of the regulari\sz ation \eqref{eq:regularized_potential}-\eqref{eq:regularized_potential2}, which is essentially the same employed in \cite{Agosti2017,Agosti2024}, with a minor  modification.
In particular, in our setting the approximation for $r<0$ of the convex part $F_1$ is specifically designed to prove a lower bound for the phase variable $\varphi$ (see Subsection \ref{subsec:Passage to the limit}).
\begin{figure}[h]
    \centering
    \includegraphics[scale=.6]{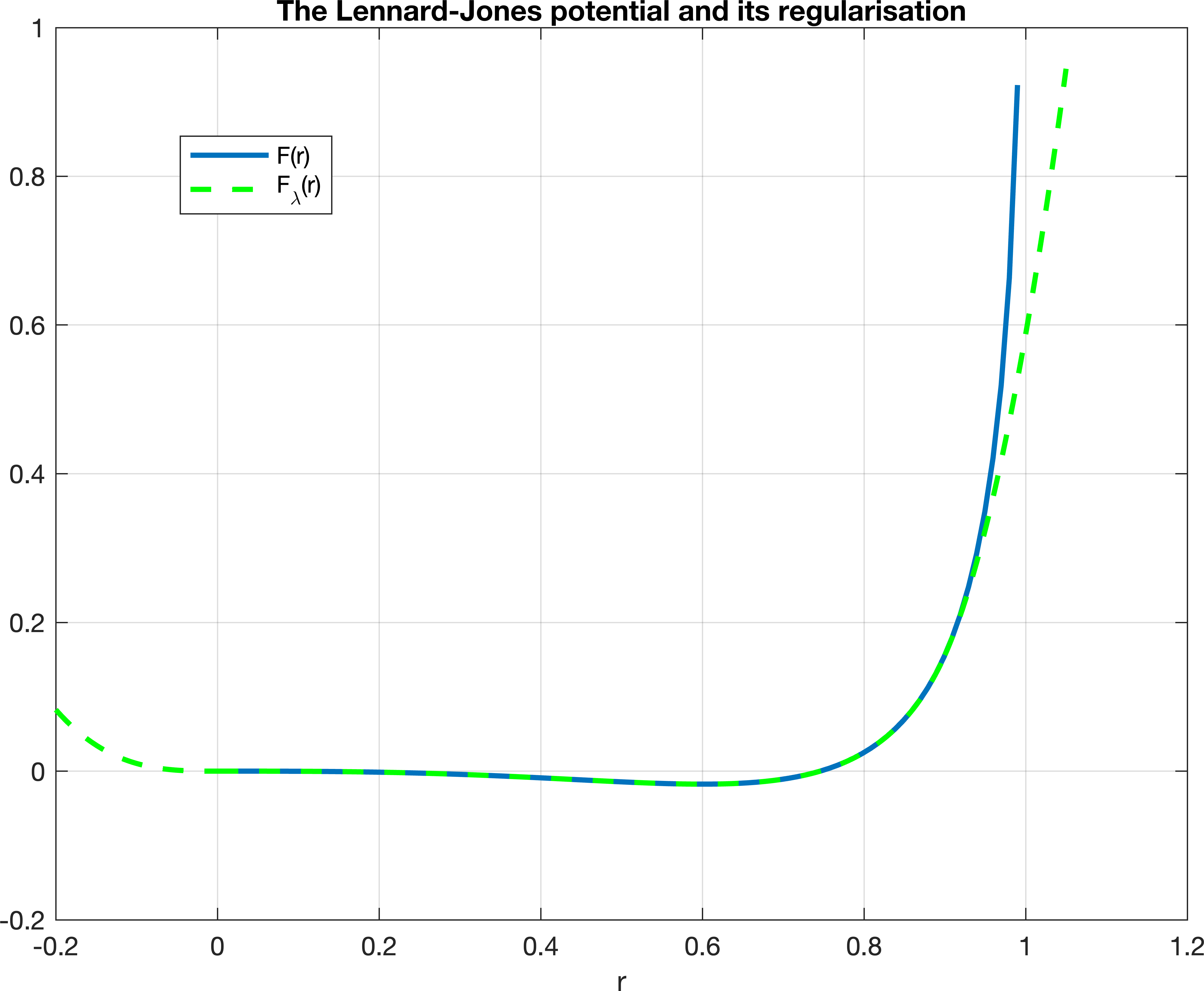} 
    \caption{The blue line represents the Lennard--Jones potential \eqref{def:LJ}, and the green dashed line its regularisation \eqref{eq:regularized_potential}-\eqref{eq:regularized_potential2}, corresponding to the value $\lambda=0.1$. Here $h$ has been selected as $h=0.4$, corresponding to $\ov\ph=0.6$.}
    \label{fig:01}
\end{figure}
From the definitions \eqref{eq:regularized_potential}-\eqref{eq:regularized_potential2}, it readily follows that $F_\lambda\to F$ and $F'_\lambda\to F'$ uniformly on all compact subsets of $[0,1)$, as $\lambda\to0^+$.
In the following, we consider $\lambda\in(0,\bar{\lambda}]$ for some $\bar{\lambda}\in(0,1)$ small enough.
For these values of $\lambda$, the approximated potential $F_\lambda$ satisfies
\begin{align*}
    |F_\lambda(r)|\geq c_1 r^2-c_2\quad\forall r\in\mathbb{R},
\end{align*}
for some constants $c_1,c_2>0$ independent of $\lambda$, where $c_1$ can be chosen arbitrarily large.
In particular, there exist two constants $\epsilon,c_2>0$, independent of $\lambda$, such that
\begin{align}
    F_\lambda(r)\geq \biggl(\frac{a^*-a_*}{2}+\epsilon\biggr) r^2-c_2 \quad \forall r\in\mathbb{R}.
    \label{eq:Flambda_bigger_than_parabola}
\end{align}
Further, we notice that
\begin{align}
    F_\lambda(r)\leq F(r)\quad\forall r\in\mathbb{R},\\
    F_\lambda''(r)+a_*\geq c_0\quad\forall r\in\mathbb{R}, \label{eq:property_of_F_second}
\end{align}
where $c_0$ is defined in \ref{ass:kernel}.
Finally, for all $\epsilon\in(0,\frac{1}{2})$, there exist constants $C_1,C_2>0$, depending on $\epsilon$ and $F$ only, such that
\begin{align}
    |F_\lambda'(r)|\leq C_1 F_\lambda'(r)(r-r_0)+C_2\quad\forall r\in\mathbb{R},\forall r_0\in[\epsilon,1-\epsilon].
    \label{eq:mz_applied_lambda}
\end{align}
To prove \eqref{eq:mz_applied_lambda}, we adapt an argument from \cite{Miranville2004} (see Proposition A.1 therein, cf. also \cite{Colli2018,Gilardi2009,Gilardi2023}).
Defining $\bar{r}:=\max\{\bar{\varphi},1-\frac{\epsilon}{2}\}$, we can divide the proof into two cases depending on $r\in[0,\bar{r}]$ or $r\in(-\infty,0)\cup(\bar{r},+\infty)$.
If $r\in[0,\bar{r}]$, then, for $\lambda$ small enough, there exists $C_2>0$, depending on $\epsilon$ and $F$ only, such that
\begin{align*}
    |F_\lambda'(r)|=|F'(r)|\leq\max_{r\in[0,\bar{r}]}|F'(r)|\leq C_2.
\end{align*}
On the other hand, if $r\in(-\infty,0)\cup(\bar{r},+\infty)$, we observe that $F_\lambda'(r)(r-r_0)>0$.
Thus we have
\begin{align*}
    |F_\lambda'(r)|=\frac{1}{|r-r_0|}F_\lambda'(r)(r-r_0)\leq C_1F_\lambda'(r)(r-r_0),
\end{align*}
for all constant $C_1$ such that $C_1\geq\frac{1}{|r-r_0|}>0$. For example, we can choose $C_1:=\frac{2}{\epsilon}$.

We can now introduce the Faedo--Galerkin approximation for our problem with $F$ replaced by $F_\lambda$.
Let $\{e_j\}_{j\in\mathbb{N}}\subset V$ and $\{l_j\}_{j\in\mathbb{N}}\subset\mathbb{R}$ be the sequences of eigenfunctions and eigenvalues of the operator $-\Delta$ with homogeneous Neumann boundary condition, respectively.
It is well known that $\{e_j\}_{j\in\mathbb{N}}$ forms an orthogonal basis for both $V$ and $H$, and we normalise it in such a way that $\|e_j\|=1$ for all $j\in\mathbb{N}$. Note also that $\{e_j\}_{j\in\mathbb{N}}\subset W$.
For each $n\in\mathbb{N}$, we set $\mathcal{W}_n:=${\rm Span}$\{e_j\}_{j=1}^n$ and denote by $\Pi_n:H\to\mathcal{W}_n$ the orthogonal projection on $\mathcal{W}_n$ with respect to the $H$-inner product $(\cdot,\cdot)$.
We then consider the following finite dimensional regularised problem: find a triplet $(\varphi_{\lambda,n},\mu_{\lambda,n},\sigma_{\lambda,n})$ 
of the form
\begin{align}
    & 
        \label{eq:Galerkin_ansatz:1} 
    \varphi_{\lambda,n}(x,t):=\sum_{j=1}^n\alpha_j^{\lambda,n}(t)e_j(x),\quad\mu_{\lambda,n}(x,t):=\sum_{j=1}^n\beta_j^{\lambda,n}(t)e_j(x),
    \\
     & \sigma_{\lambda,n}(x,t):=\sum_{j=1}^n\gamma_j^{\lambda,n}(t)e_j(x),
    \label{eq:Galerkin_ansatz:2}
\end{align}
where $t\in[0,T]$ and $x\in\Omega$, such that
\begin{alignat}{2}
    &\partial_t\varphi_{\lambda,n}-\Delta\mu_{\lambda,n}=\Pi_n[S(\varphi_{\lambda,n},\sigma_{\lambda,n})]&&\quad \text{ in }Q, \label{eq:varphi_lambda_n} \\
    &\mu_{\lambda,n}=\tau\partial_t\varphi_{\lambda,n}+a\varphi_{\lambda,n}-J\ast\varphi_{\lambda,n}+\Pi_nF_\lambda'(\varphi_{\lambda,n})-\chi\sigma_{\lambda,n}&&\quad \text{ in }Q, \label{eq:mu_lambda_n} \\
    &\partial_t\sigma_{\lambda,n}-\Delta\sigma_{\lambda,n}+B(\sigma_{\lambda,n}-\Pi_n\sigma_{S})+\Pi_n[C\sigma_{\lambda,n} {\hhd}(\varphi_{\lambda,n})]=0&&\quad \text{ in }Q, \label{eq:sigma_lambda_n} \\
    &\partial_{\boldsymbol{n}}\mu_{\lambda,n}=\partial_{\boldsymbol{n}}\sigma_{\lambda,n}=0&&\quad \text{ on }\Sigma, \label{eq:bcs_lambda_n} \\
    &\varphi_{\lambda,n}(0)=\Pi_n\varphi_0,\quad\sigma_{\lambda,n}(0)=\Pi_n\sigma_0&&\quad \text{ in }\Omega. \label{eq:ics_lambda_n}
\end{alignat}
Thus, the problem is to find
$\alpha_j^{\lambda,n},\beta_j^{\lambda,n},\gamma_j^{\lambda,n}:[0,T]\to\mathbb{R}$, for all $j\in\{1,...,n\}$ and for all $(\lambda,n)\in(0,1)\times\mathbb{N}$.
It is then natural to introduce the vector-valued functions
\begin{align*}
    \boldsymbol{\alpha}^{\lambda,n},\boldsymbol{\beta}^{\lambda,n},\boldsymbol{\gamma}^{\lambda,n}:[0,T]\to\mathbb{R}^n,
\end{align*}
where
\begin{align*}
    (\boldsymbol{\alpha}^{\lambda,n})_j=\alpha_j^{\lambda,n},\quad(\boldsymbol{\beta}^{\lambda,n})_j=\beta_j^{\lambda,n},\quad(\boldsymbol{\gamma}^{\lambda,n})_j=\gamma_j^{\lambda,n}\qquad\forall j\in\{1,...,n\}.
\end{align*}
By a standard procedure, inserting \eqref{eq:Galerkin_ansatz:1}-\eqref{eq:Galerkin_ansatz:2} into \eqref{eq:varphi_lambda_n}-\eqref{eq:ics_lambda_n} and testing by $e_i$ for each $i\in\{1,...,n\}$, we see that the triplet $(\varphi_{\lambda,n},\mu_{\lambda,n},\sigma_{\lambda,n})$ solves \eqref{eq:varphi_lambda_n}-\eqref{eq:ics_lambda_n} if and only if the vectors $\boldsymbol{\alpha}^{\lambda,n},\boldsymbol{\beta}^{\lambda,n},\boldsymbol{\gamma}^{\lambda,n}$ solve a Cauchy problem of the form
\begin{equation*}
    \begin{cases}
        \diff{}{t}(\boldsymbol{\alpha}^{\lambda,n},\boldsymbol{\beta}^{\lambda,n},\boldsymbol{\gamma}^{\lambda,n}) =\boldsymbol{g}_{\lambda,n}(\boldsymbol{\alpha}^{\lambda,n},\boldsymbol{\beta}^{\lambda,n},\boldsymbol{\gamma}^{\lambda,n}),\\
        (\boldsymbol{\alpha}^{\lambda,n},\boldsymbol{\beta}^{\lambda,n},\boldsymbol{\gamma}^{\lambda,n})(0)=(\{(\varphi_0,e_i)\}_{i=1}^n,\{(\mu_0,e_i)\}_{i=1}^n,\{(\sigma_0,e_i)\}_{i=1}^n),
    \end{cases}
\end{equation*}
where $\boldsymbol{g}_{\lambda,n}:\mathbb{R}^{3n}\to\mathbb{R}^{3n}$ is continuous in view of the continuity of ${\hhu},{\hhd}$ and $F_\lambda'$.
Thus, there exists a local solution $(\boldsymbol{\alpha}^{\lambda,n},\boldsymbol{\beta}^{\lambda,n},\boldsymbol{\gamma}^{\lambda,n})\in\mathcal{C}^1([0,T^*];\mathbb{R}^{3n})$ for some maximal time  $T^*\in(0,T]$. Recalling \eqref{eq:Galerkin_ansatz:1}-\eqref{eq:Galerkin_ansatz:2}, this entails the existence of $(\varphi_{\lambda,n},\mu_{\lambda,n},\sigma_{\lambda,n})\in\mathcal{C}^1([0,T^*];\mathcal{W}_n)$ that solves \eqref{eq:varphi_lambda_n}-\eqref{eq:ics_lambda_n} for all $t\in[0,T^*]$.


\subsection{A priori estimates}
\label{subsec:Uniform estimates}
Here we prove some a priori estimates that are uniform with respect to both $n$ and $\lambda$. These estimates will allow us to show that $T^*=T$, that is that the local solution $(\varphi_{\lambda,n},\mu_{\lambda,n},\sigma_{\lambda,n})$ can be extended to the whole time interval $[0,T]$. Furthermore, they will enable us to to extract a suitable subsequence which properly converges to a solution to problem \Sys.

We multiply \eqref{eq:varphi_lambda_n} by $\mu_{\lambda,n}$, \eqref{eq:mu_lambda_n} by $-\partial_t\varphi_{\lambda,n}$ and \eqref{eq:sigma_lambda_n} by $\sigma_{\lambda,n}$, take the sum and integrate over \modLuca{$Q_t$} for any $t\in[0,T]$.
This gives
\begin{align*}
    & \int_{Q_t}|\nabla\mu_{\lambda,n}|^2+\tau\int_{Q_t}|\partial_t\varphi_{\lambda,n}|^2 
+\frac{1}{4}\int_{\Omega\times\Omega}J(x-y)|\varphi_{\lambda,n}(x,t)-\varphi_{\lambda,n}(y,t)|^2\dd x\dd y
    \\
    &\qquad 
    +\int_\Omega F_\lambda(\varphi_{\lambda,n}(t))+\frac{1}{2}\|\sigma_{\lambda,n}(t)\|^2+\int_{Q_t}|\nabla\sigma_{\lambda,n}|^2+\int_{Q_t}(B+C{\hhd}(\varphi_{\lambda,n}))|\sigma_{\lambda,n}|^2\\
    & \quad = \frac{1}{4}\int_{\Omega\times\Omega}J(x-y)|\Pi_n\varphi_0(x)-\Pi_n\varphi_0(y)|^2\dd x\dd y+\int_\Omega F_\lambda(\Pi_n\varphi_0)\\
    &\qquad +\frac{1}{2}\|\Pi_n\sigma_0\|^2+\int_{Q_t}S(\varphi_{\lambda,n},\sigma_{\lambda,n})\mu_{\lambda,n}+\chi\int_{Q_t}\sigma_{\lambda,n}\partial_t\varphi_{\lambda,n}+B\int_{Q_t}\sigma_{S}\sigma_{\lambda,n}.
\end{align*}
Using now the Young inequality, the Young convolution inequality, \eqref{eq:Flambda_bigger_than_parabola}, and recalling \ref{eq:assumption3} and \ref{ass:kernel}, we find
\begin{align}
	\non
	    & \int_{Q_t}|\nabla\mu_{\lambda,n}|^2+\tau\int_{Q_t}|\partial_t\varphi_{\lambda,n}|^2+\epsilon\|\varphi_{\lambda,n}(t)\|^2+\frac{1}{2}\|\sigma_{\lambda,n}(t)\|^2+\int_{Q_t}|\nabla\sigma_{\lambda,n}|^2\\
   &
   \non \quad  \leq  a^*\|\varphi_0\|^2+\frac{1}{2}\|\sigma_{0}\|^2+\int_\Omega F_\lambda(\Pi_n\varphi_0)+M\biggl(1+\int_{Q_t}|\sigma_{\lambda,n}|^2\biggr)\\
    &\qquad +\chi\int_{Q_t}\sigma_{\lambda,n}\partial_t\varphi_{\lambda,n}+\int_{Q_t}S(\varphi_{\lambda,n},\sigma_{\lambda,n})\mu_{\lambda,n},
    \label{eq:53_0}
\end{align}
for some $\epsilon>0$.
Next, let us multiply \eqref{eq:varphi_lambda_n} by $\varphi_{\lambda,n}$ and \eqref{eq:mu_lambda_n} by $\Delta\varphi_{\lambda,n}$. Adding the resulting equalities and then integrating over \modLuca{$Q_t$}, we get
\begin{align}
\begin{split}
    \frac{1}{2}\|\varphi&_{\lambda,n}(t)\|^2+\frac{\tau}{2}\|\nabla\varphi_{\lambda,n}(t)\|^2\\
    =&\int_{Q_t}S(\varphi_{\lambda,n},\sigma_{\lambda,n})\varphi_{\lambda,n}+\chi\int_{Q_t}\nabla\sigma_{\lambda,n}\cdot\nabla\varphi_{\lambda,n}+\int_{Q_t}F_\lambda'(\varphi_{\lambda,n})\Delta\varphi_{\lambda,n}\\
    &+\int_{Q_t}(a\varphi_{\lambda,n}-J\ast\varphi_{\lambda,n})\Delta\varphi_{\lambda,n}+\frac{1}{2}\|\Pi_n\varphi_0\|^2+\frac{\tau}{2}\|\nabla\Pi_n\varphi_0\|^2.
    \label{eq:second_estimate_equality_start}
\end{split}
\end{align}
The first term on the right-hand side of \eqref{eq:second_estimate_equality_start} can be estimated using the Young inequality and taking \ref{eq:assumption1} into account. The second one can be estimated by means of the Young inequality as follows:
\begin{align}
    \chi\int_{Q_t}\nabla\varphi_{\lambda,n}\cdot\nabla\sigma_{\lambda,n}\leq\frac{\chi^2}{2}\int_{Q_t}|\nabla\varphi_{\lambda,n}|^2+\frac{1}{2}\int_{Q_t}|\nabla\sigma_{\lambda,n}|^2.
    \label{eq:grad_grad_estimate}
\end{align}
Then, we can integrate by parts the third term on the right-hand side of \eqref{eq:second_estimate_equality_start} and use \eqref{eq:property_of_F_second}:
\begin{align*}
    \int_{Q_t}F_\lambda'(\varphi_{\lambda,n})\Delta\varphi_{\lambda,n}=-\int_{Q_t}F_\lambda''(\varphi_{\lambda,n})|\nabla\varphi_{\lambda,n}|^2\leq-\int_{Q_t}(c_0-a_*)|\nabla\varphi_{\lambda,n}|^2.
\end{align*}
Finally, the fourth term on the right-hand side of \eqref{eq:second_estimate_equality_start} can be treated by means of an integration by parts, the Hölder inequality, assumption \ref{ass:kernel}, the Young convolution inequality and the Young inequality. The resulting inequality is
\begin{align*}
    \int_{Q_t}(a\varphi_{\lambda,n}-J\ast\varphi_{\lambda,n})\Delta\varphi_{\lambda,n}\leq\biggl(\frac{c_0}{2}-a^*\biggr)\int_{Q_t}|\nabla\varphi_{\lambda,n}|^2+\frac{b^*}{2c_0}\int_{Q_t}|\varphi_{\lambda,n}|^2.
\end{align*}
Taking into account the inequalities above and exploiting the continuity of $\Pi_n$, the fact that $a^*-a_*\geq0$, and the triangular inequality, we infer from \eqref{eq:second_estimate_equality_start} that
\begin{align}
\begin{split}
    \frac{1}{2}\|\varphi&_{\lambda,n}(t)\|^2+\frac{\tau}{2}\|\nabla\varphi_{\lambda,n}(t)\|^2+\frac{c_0}{2}\int_{Q_t}|\nabla\varphi_{\lambda,n}|^2\\
    \leq& M\biggl(1+\int_{Q_t}|\varphi_{\lambda,n}|^2\biggr)+\frac{1}{2}\int_{Q_t}|\nabla\sigma_{\lambda,n}|^2+\|\varphi_0\|^2+\frac{\tau}{2}\|\varphi_0\|_V^2+\frac{\chi^2}{2}\int_{Q_t}|\nabla\varphi_{\lambda,n}|^2.
    \label{eq:53_1}
\end{split}
\end{align}
Summing \eqref{eq:53_0} to \eqref{eq:53_1} yields
\begin{align}
\begin{split}
    \int_{Q_t}|\nabla\mu&_{\lambda,n}|^2+\tau\int_{Q_t}|\partial_t\varphi_{\lambda,n}|^2+\epsilon\|\varphi_{\lambda,n}(t)\|^2+\frac{1}{2}\|\sigma_{\lambda,n}(t)\|^2+\frac{1}{2}\int_{Q_t}|\nabla\sigma_{\lambda,n}|^2\\
    &+\frac{1}{2}\|\varphi_{\lambda,n}(t)\|^2+\frac{\tau}{2}\|\nabla\varphi_{\lambda,n}(t)\|^2+\frac{c_0}{2}\int_{Q_t}|\nabla\varphi_{\lambda,n}|^2\\
    \leq& a^*\|\varphi_0\|^2+\frac{1}{2}\|\sigma_{0}\|^2+\|\varphi_{0}\|^2+\frac{\tau}{2}\|\varphi_{0}\|_V^2+\int_\Omega F_\lambda(\Pi_n\varphi_0)+\frac{\chi^2}{2}\int_{Q_t}|\nabla\varphi_{\lambda,n}|\\
    &+M\biggl(1+\int_{Q_t}|\varphi_{\lambda,n}|^2+\int_{Q_t}|\sigma_{\lambda,n}|^2\biggr)+\chi\int_{Q_t}\sigma_{\lambda,n}\partial_t\varphi_{\lambda,n}+\int_{Q_t}S(\varphi_{\lambda,n},\sigma_{\lambda,n})\mu_{\lambda,n}.
    \label{eq:fondamentale}
\end{split}
\end{align}
The last term on the right-hand side of \eqref{eq:fondamentale} can be estimated recalling equation \eqref{eq:mu_lambda_n} and the fact that $\tau\leq1$, using the Young inequality. This gives
\begin{align}
\begin{split}
    \int_{Q_t}S(\varphi&_{\lambda,n},\sigma_{\lambda,n})\mu_{\lambda,n}\\
    =&\tau\int_{Q_t} S(\varphi_{\lambda,n},\sigma_{\lambda,n})\partial_t\varphi_{\lambda,n}+a\int_{Q_t} S(\varphi_{\lambda,n},\sigma_{\lambda,n})\varphi_{\lambda,n}-\int_{Q_t} S(\varphi_{\lambda,n},\sigma_{\lambda,n})J\ast\varphi_{\lambda,n}\\
    &+\int_{Q_t}S(\varphi_{\lambda,n},\sigma_{\lambda,n})F_{1,\lambda}'(\varphi_{\lambda,n})+\int_{Q_t}S(\varphi_{\lambda,n},\sigma_{\lambda,n})\bar{F}_2'(\varphi_{\lambda,n})-\chi\int_{Q_t}S(\varphi_{\lambda,n},\sigma_{\lambda,n})\sigma_{\lambda,n}\\
    \leq&\int_{Q_t} S(\varphi_{\lambda,n},\sigma_{\lambda,n})F_{1,\lambda}'(\varphi_{\lambda,n})-m\int_{Q_t}\varphi_{\lambda,n}\bar{F}_2'(\varphi_{\lambda,n})+\int_{Q_t}{\hhu}(\varphi_{\lambda,n},\sigma_{\lambda,n})\bar{F}_2'(\varphi_{\lambda,n})\\
    &+\frac{\tau}{2}\int_{Q_t}|\partial_t\varphi_{\lambda,n}|^2+M\biggl(1+\int_{Q_t}|\varphi_{\lambda,n}|^2+\int_{Q_t}|\sigma_{\lambda,n}|^2\biggr).
    \label{eq:fondamentale2}
\end{split}
\end{align}
Recalling that $F_{1,\lambda}$ is a convex function, the first term on the right-hand side of  \eqref{eq:fondamentale2} can be bounded by a constant as follows
\begin{align*}
    \int_{Q_t}&S(\varphi_{\lambda,n},\sigma_{\lambda,n})F_{1,\lambda}'(\varphi_{\lambda,n})=m\int_{Q_t}\biggl(\frac{{\hhu}(\varphi_{\lambda,n},\sigma_{\lambda,n})}{m}-\varphi_{\lambda,n}\biggr)F_{1,\lambda}'(\varphi_{\lambda,n})\\
    &\leq m\int_{Q_t}F_{1,\lambda}\biggl(\frac{{\hhu}(\varphi_{\lambda,n},\sigma_{\lambda,n})}{m}\biggr)-m\int_{Q_t}F_{1,\lambda}(\varphi_{\lambda,n})\leq M,
\end{align*}
where we have used (see \ref{eq:assumption1}) 
\begin{align*}
    \frac{{\hhu}(\varphi_{\lambda,n},\sigma_{\lambda,n})}{m}\in \bigl[0,\tfrac{K}{m}\bigr)\subset[0,1)\quad\text{a.e. in }Q,
\end{align*}
the fact that $\lambda$ can be chosen arbitrarily close to zero and that $F_{1,\lambda}$ is bounded below independently of $\lambda$.
Regarding the second term on the right-hand side of inequality \eqref{eq:fondamentale2}, we observe that the function $r\mapsto-r\bar{F}_2'(r)$ is negative if $r<0$, is continuous for all $r\in\mathbb{R}$, and has quadratic growth for $r>1$.
These three facts readily imply that there exists $M>0$ such that $-r\bar{F}_2'(r)\leq M(1+|r|^2)$ for all $r\in\mathbb{R}$. Therefore, we have
\begin{align*}
    -m\int_{Q_t}\varphi_{\lambda,n}\bar{F}_2'(\varphi_{\lambda,n})\leq M\biggl(1+\int_{Q_t}|\varphi_{\lambda,n}|^2\biggr).
\end{align*}
Finally, exploiting \ref{eq:assumption1} and the fact that $\bar{F}_2'$ has at most a quadratic growth, we get the existence of a constant $M>0$ that allows to bound the third term on the right-hand side of inequality \eqref{eq:fondamentale2}:
\begin{align*}
    \int_{Q_t}{\hhu}(\varphi_{\lambda,n},\sigma_{\lambda,n})\bar{F}_2'(\varphi_{\lambda,n})\leq M\biggl(1+\int_{Q_t}|\varphi_{\lambda,n}|^2\biggr).
\end{align*}
Applying the Young inequality to the eighth term on the right-hand side of \eqref{eq:fondamentale}, taking into account all the previous remarks, and using \eqref{eq:fondamentale2}, we get from \eqref{eq:fondamentale} that
\begin{align}
\begin{split}
    \int_{Q_t}|\nabla\mu&_{\lambda,n}|^2+\frac{\tau}{4}\int_{Q_t}|\partial_t\varphi_{\lambda,n}|^2+\epsilon\|\varphi_{\lambda,n}(t)\|^2+\frac{1}{2}\|\sigma_{\lambda,n}(t)\|^2+\frac{1}{2}\int_{Q_t}|\nabla\sigma_{\lambda,n}|^2\\
    &+\frac{1}{2}\|\varphi_{\lambda,n}(t)\|^2+\frac{\tau}{2}\|\nabla\varphi_{\lambda,n}(t)\|^2+\frac{c_0}{2}\int_{Q_t}|\nabla\varphi_{\lambda,n}|^2\\
    \leq&M_{\tau}\biggl(1+\int_{Q_t}|\varphi_{\lambda,n}|^2+\int_{Q_t}|\sigma_{\lambda,n}|^2+\int_{Q_t}|\nabla\varphi_{\lambda,n}|\biggr)+\int_\Omega F_\lambda(\Pi_n\varphi_0),
    \label{eq:53_4}
\end{split}
\end{align}
in view of \eqref{eq:assumptions_ics}.
Noticing that the last term on the right-hand side of \eqref{eq:53_4} is bounded independently of $n$, the Gronwall lemma yields
\begin{align}
    \|\varphi_{\lambda,n}\|_{H^1(0,T;H)\cap L^\infty(0,T;V)}+\|\nabla\mu_{\lambda,n}\|_{L^2(0,T;H)}+\|\sigma_{\lambda,n}\|_{L^\infty(0,T;H)\cap L^2(0,T;V)}\leq M_{\lambda,\tau}.
    \label{eq:53_5}
\end{align}
In order to bound the $L^2(0,T;H)$-norm of $\mu_{\lambda,n}$,
we \modLuca{test} \eqref{eq:mu_lambda_n} with $|\Omega|^{-1}$. Then, taking the $L^2(0,T)$-norm and exploiting the triangular inequality and \eqref{eq:53_5}, we obtain
\begin{align*}
    \|(\mu_{\lambda,n})_\Omega\|_{L^2(0,T)}\leq M_{\lambda,\tau}\biggl(1+\|F'_\lambda(\varphi_{\lambda,n})\|_{L^2(0,T;L^1(\Omega))}\biggr).
\end{align*}
Since $F_\lambda'$ has sub-quadratic growth, owing to \eqref{eq:53_5}, we have that
\begin{align*}
    \|F'_\lambda(\varphi_{\lambda,n})\|_{L^2(0,T;L^1(\Omega))}\leq M_\lambda\biggl(1+\|\varphi_{\lambda,n}\|^2_{L^4(0,T;H)}\biggr)\leq M_{\lambda,\tau}.
\end{align*}
Hence, we have $\|(\mu_{\lambda,n})_\Omega\|_{L^2(0,T)}\leq M_{\lambda,\tau}$ so that, owing the Poincar\'e--Wirtinger inequality, we obtain
\begin{align*}
    \|\mu_{\lambda,n}\|_{L^2(0,T;V)}\leq M_{\lambda,\tau}.
\end{align*}
Finally, by comparison in \eqref{eq:sigma_lambda_n}, recalling \ref{eq:assumption2}-\ref{eq:assumption3}, and using \eqref{eq:53_5}, we infer
\begin{align*}
    \|\sigma_{\lambda,n}\|_{H^1(0,T;V')}\leq M_{\lambda,\tau}.
\end{align*}

\subsection{Passage to the limit}
\label{subsec:Passage to the limit}

We now pass to the limit, first as $n \to \infty$ and then as $\lambda \to 0$, along suitable subsequences. We collect the above estimates to obtain
\begin{align*}
     & \|\varphi_{\lambda,n}\|_{H^1(0,T;H) \cap L^\infty(0,T;V)} + \|\mu_{\lambda,n}\|_{L^2(0,T;V)} 
     + \|\sigma_{\lambda,n}\|_{H^1(0,T;V') \cap L^\infty(0,T;H) \cap L^2(0,T;V)} \leq M_{\lambda,\tau}.
\end{align*}
Thus, using standard compactness arguments and the Aubin--Lions compactness results (see, e.g., \cite[Cor. 4]{simon}), we deduce that there exists a limiting
triplet $(\varphi_\lambda,\mu_\lambda,\sigma_\lambda)$ such that, up to subsequences, as $n \to \infty$,
\begin{align*}
    \ph_{\lambda,n} \overset{*}{\rightharpoonup} \varphi_\lambda&\qquad  \text{in $H^1(0,T;H)\cap L^\infty(0,T;V),$}
    \quad 
    \ph_{\lambda,n} \to  \varphi_{\lambda}\quad  \text{in }\mathcal{C}^0([0,T];H),\\
    \mu_{\lambda,n} {\rightharpoonup} \mu_\lambda &\qquad \text{in $L^2(0,T;V)$},\\
    \sigma_{\lambda,n} \overset{*}{\rightharpoonup} \sigma_\lambda& \qquad \text{in  $H^1(0,T;V')\cap\L\infty H\cap L^2(0,T;V)$},
    \quad \sigma_{\lambda,n} \to \sigma_{\lambda} \quad \text{in  $\L2 H$},
\end{align*}
and $(\varphi_\lambda,\mu_\lambda,\sigma_\lambda)$ solves

    \begin{align}
       & \langle \partial_t\varphi_\lambda,\xi \rangle + \int_\Omega\nabla\mu_\lambda\cdot\nabla\xi=\int_\Omega S(\varphi_\lambda,\sigma_\lambda)\xi, \label{eq:53_7} \\
        &\langle \partial_t\sigma_\lambda,\xi \rangle + \int_\Omega\nabla\sigma_\lambda\cdot\nabla\xi+\int_\Omega (B(\sigma_\lambda-\sigma_S)+C\sigma_\lambda {\hhd}(\varphi_\lambda))\xi=0, 
    \end{align}
    for all $\xi \in V$, almost everywhere in $(0,T)$, with 
    \begin{alignat}{2}
        & \mu_\lambda=\tau\partial_t\varphi_\lambda+a\varphi_\lambda-J\ast\varphi_\lambda+F_\lambda'(\varphi_\lambda)-\chi\sigma_\lambda &&\qquad\text{a.e. in }Q, \label{eq:mueps_lambda} \\
        & \varphi_\lambda(0)=\varphi_0,\quad\sigma_\lambda(0)=\sigma_0 &&\qquad\text{a.e. in }\Omega. \label{eq:53_75}
    \end{alignat}

Recalling Subsection \ref{subsec:a_priori_estimates}, we can prove that there is $\delta\in(0,\frac{1}{2})$, independent of $\lambda$, such that
\begin{align}
    \delta\leq(\varphi_{\lambda})_\Omega(t)\leq1-\delta\quad\forall t\in[0,T].
    \label{eq:53_8}
\end{align}
In this setting, we can retrieve a bound similar to \eqref{eq:53_4}, namely,
\begin{align}
    \|\varphi_{\lambda}\|_{H^1(0,T;H)\cap L^\infty(0,T;V)}+\|\nabla\mu_{\lambda}\|_{L^2(0,T;H)}+\|\sigma_{\lambda}\|_{L^\infty(0,T;H)\cap L^2(0,T;V)}\leq M_{\tau}.
    \label{eq:53_9}
\end{align}

To control the $L^2(0,T;H)$-norm of the chemical potential independently of $\lambda$, we can integrate \eqref{eq:mueps_lambda} over $\Omega$ to get
\begin{align*}
    |(\mu_\lambda)_\Omega(s)|\leq M_\tau\biggl(\|\partial_t\varphi_\lambda(s)\|_{L^1(\Omega)}+\|F_\lambda'(\varphi_\lambda(s))\|_{L^1(\Omega)}+\|\sigma_\lambda(s)\|_{L^1(\Omega)}\biggr),
\end{align*}
for almost every $s\in(0,T)$.
Then, taking \eqref{eq:53_9} into account, we infer
\begin{align*}
    \|(\mu_\lambda)_\Omega\|^2_{L^2(0,T)}\leq M_\tau\biggl(1+\|F_\lambda'(\varphi_\lambda)\|^2_{L^2(0,T;L^1(\Omega))}\biggr).
\end{align*}
To control the term on the right-hand side, we notice that, recalling \eqref{eq:mz_applied_lambda} and \eqref{eq:53_8}, we can find $C_1,C_2>0$ such that
\begin{align*}
    \|F_\lambda'(\varphi_\lambda)\|_{L^1(\Omega)}\leq C_1\int_\Omega F_\lambda'(\varphi_\lambda)(\varphi_\lambda-(\varphi_\lambda)_\Omega)+C_2|\Omega|,
\end{align*}
almost everywhere in $(0,T)$.
Moreover, using \eqref{eq:mueps_lambda} and the identity $\int_{\Omega}(\mu_\lambda)_\Omega(\varphi_\lambda-(\varphi_\lambda)_\Omega)=0$, we arrive at
\begin{align}
\begin{split}
    \|F_\lambda'(\varphi_\lambda)\|_{L^1(\Omega)}\leq M\biggl(&\int_\Omega (\mu_\lambda-(\mu_\lambda)_\Omega)(\varphi_\lambda-(\varphi_\lambda)_\Omega)-\tau\int_\Omega\partial_t\varphi_\lambda(\varphi_\lambda-(\varphi_\lambda)_\Omega) \\
    &-\int_\Omega(a\varphi_\lambda-J*\varphi_\lambda)(\varphi_\lambda-(\varphi_\lambda)_\Omega)+\int_\Omega\sigma_\lambda(\varphi_\lambda-(\varphi_\lambda)_\Omega)\biggr)+M.
    \label{eq:53_10}
\end{split}
\end{align}
The first term on the right-hand side of \eqref{eq:53_10} can be bounded as follows
\begin{align}
    \int_\Omega (\mu_\lambda-(\mu_\lambda)_\Omega)(\varphi_\lambda-(\varphi_\lambda)_\Omega)\leq\|\mu_\lambda-(\mu_\lambda)_\Omega\|\|\varphi_\lambda-(\varphi_\lambda)_\Omega\|\leq M\|\nabla\mu_\lambda\|\|\varphi_\lambda\|.
    \label{eq:application_of_PW}
\end{align}
On the other hand, we have
\begin{align*}
    -\tau \int_\Omega\partial_t\varphi_\lambda(\varphi_\lambda-(\varphi_\lambda)_\Omega)\leq M_\tau\|\partial_t\varphi_\lambda\|(1+\|\varphi_\lambda\|).
\end{align*}
Then, thanks to \eqref{eq:53_8}, we have
\begin{align*}
    -\int_\Omega(a\varphi_\lambda-J*\varphi_\lambda)(\varphi_\lambda-(\varphi_\lambda)_\Omega)\leq M(1+\|\varphi_\lambda\|^2),
\end{align*}
and
\begin{align*}
    \int_\Omega\sigma_\lambda(\varphi_\lambda-(\varphi_\lambda)_\Omega)\leq M(1+\|\varphi_\lambda\|^2+\|\sigma_\lambda\|^2).
\end{align*}
Taking these facts into account, from \eqref{eq:53_10} we derive that
\begin{align*}
    \|F_\lambda'(\varphi_\lambda)\|_{L^1(\Omega)}\leq M_\tau\bigl(1+\|\nabla\mu_\lambda\|\|\varphi_\lambda\|+\|\partial_t\varphi_\lambda\|(1+\|\varphi_\lambda\|)+\|\sigma_\lambda\|^2+\|\varphi_\lambda\|^2\bigr),
\end{align*}
which gives
\begin{align*}
    \modLuca{\|F_\lambda'(\varphi_\lambda)\|^2_{\L2 {\Lx1}}}\leq M_\tau\bigl(1&+\|\varphi_\lambda\|^2_{L^\infty(0,T;H)}\|\nabla\mu_\lambda\|^2_{L^2(0,T;H)}+\|\varphi_\lambda\|^2_{L^\infty(0,T;H)}\\
    &+\|\partial_t\varphi_\lambda\|^2_{L^2(0,T;H)}\bigl(1+\|\varphi_\lambda\|^2_{L^\infty(0,T;H)}\bigr)+\|\sigma_\lambda\|^2_{L^\infty(0,T;H)}\bigr),
\end{align*}
where the right-hand side is bounded independently of $\lambda$ in view of \eqref{eq:53_9}.
This allows to conclude that
\begin{align*}
    \|(\mu_\lambda)_\Omega\|_{L^2(0,T)}\leq M_\tau,
\end{align*}
yielding the desired $\lambda$-independent bound for $\mu_\lambda$ in the $L^2(0,T;V)$\modLuca{-norm}.
Gathering the above estimates, we obtain 
\begin{align*}
     & \|\varphi_{\lambda}\|_{H^1(0,T;H) \cap L^\infty(0,T;V)} + \|\mu_{\lambda}\|_{L^2(0,T;V)}
     + \|\sigma_{\lambda}\|_{H^1(0,T;V') \cap \C0 H \cap L^2(0,T;V)} \leq M_{\tau}.
\end{align*}
Then, arguing as above via compactness,
we find a
triplet $(\varphi,\mu,\sigma)$ such that, as $\lambda\to 0$ along a suitable subsequence,
\begin{align*}
    \ph_\lambda \overset{*}{\rightharpoonup}  \varphi\qquad  & \text{in $H^1(0,T;H)\cap L^\infty(0,T;V),$}\quad 
    \ph_\lambda \to  \varphi\quad  \text{in }\mathcal{C}^0([0,T];H),\\
    \mu_\lambda {\rightharpoonup} \mu \qquad & \text{in $L^2(0,T;V)$},\\
    \sigma_\lambda \overset{*}{\rightharpoonup} \sigma \qquad & \text{in  $H^1(0,T;V')\cap\L\infty H\cap L^2(0,T;V)$},
    \quad
    \sigma_\lambda \to \sigma \quad \text{in  $\L2 H$}.
\end{align*}
To prove that such a triplet is a weak solution to \eqref{eq:varphi}-\eqref{eq:ics} in the sense of Definition \ref{def:weak_sol}, one has to let $\lambda\to0^+$ in \eqref{eq:53_7}-\eqref{eq:53_75} along a suitable subsequence, exploiting the above convergences.
The only delicate step of this argument, 
on account of \cite[Lemma 1.3]{Lions1969}, is to show that
\begin{align}
    F_\lambda'(\varphi_\lambda)\to F'(\varphi)\quad\text{a.e. in }Q.
    \label{eq:passaggio_al_limite_critico}
\end{align}
Observe that $\varphi_{\lambda}\to\varphi\text{ in }\mathcal{C}^0([0,T];H)$, up to a subsequence, by the Aubin--Lions Lemma. Thus, it also converges almost everywhere in $Q$, up to a subsequence.
Consequently, if we claim that \eqref{eq:MP_varphi} holds, then the pointwise convergence of $\varphi_\lambda$ along with the uniform convergence of $F'_\lambda$ to $F'$ in every compact set in $[0,1)$, entails \eqref{eq:passaggio_al_limite_critico}.

To conclude the proof of Theorem \ref{thm:existence_weak_tau>0}, it suffices to prove \eqref{eq:MP_varphi}.

First, we show that $\varphi<1$ almost everywhere in $Q$.
Indeed, we can adapt Step $2$ of \cite[Subsection 3.3] {Giorgini2018} to the case of the single-well potential \eqref{def:LJ}.
Let $\Tilde{\varphi}\in(0,1)$ be the real number such that $F(\Tilde{\varphi})=0$, then fix $\rho\in(0,1-\Tilde{\varphi})$ and define the following sets:
\begin{align*}
    E_\rho^\lambda:=\bigl\{(x,t)\in Q\mid\varphi_\lambda(x,t)>1-\rho\bigr\},\qquad E_\rho:=\bigl\{(x,t)\in Q\mid\varphi(x,t)>1-\rho\bigr\}.
\end{align*}
Since $\varphi_\lambda$ converges pointwise to $\varphi$, as $\lambda\to 0$, the Fatou Lemma yields
\begin{align}
    |E_\rho|\leq\liminf_{\lambda\to0^+}{|E_\rho^\lambda|}.
    \label{eq:GGW-Fatou}
\end{align}
Moreover, taking the $L^1(Q)$-norm of \eqref{eq:mueps_lambda} and using the uniform bounds above, we obtain
\begin{align*}
    \int_Q|F_{\lambda}'(\varphi_\lambda)|\leq M_{\tau}.
\end{align*}
It then follows that
\begin{align*}
    M_{\tau}\geq \int_Q|F_{\lambda}'(\varphi_\lambda)|\geq\int_{E_\rho^\lambda}|F_{\lambda}'(\varphi_\lambda)|\geq\int_{E_\rho^\lambda}\inf_{(x,t)\in E_\rho^\lambda}|F_{\lambda}'(\varphi_\lambda(x,t))|=|E_\rho^\lambda|\min_{r\geq1-\rho}|F_{\lambda}'(r)|.
\end{align*}
We now exploit the fact that $\lambda$ can be taken as close to zero as needed, in such a way that $\lambda\leq\rho$, and the fact that $F_\lambda\in\mathcal{C}^2(\mathbb{R})$ is increasing and positive in $(1-\rho,+\infty)$, where $\rho\in(0,1-\Tilde{\varphi})$. 
Hence, we have
\begin{align*}
    \min_{r\geq1-\rho}|F_{\lambda}'(r)| = F_{\lambda}'(1-\rho) = F'(1-\rho),
\end{align*}
so that
\begin{align*}
    |E_\rho^\lambda|\leq\frac{M_{\tau}}{F'(1-\rho)}.
\end{align*}
Now, we first let $\lambda\to0^+$ using \eqref{eq:GGW-Fatou} and then $\rho\to0^+$ to deduce that
\begin{align*}
    \biggl|\biggl\{(x,t)\in Q\mid\varphi(x,t)\geq1\biggr\}\biggr|=0,
\end{align*}
i.e., $\varphi<1$ almost everywhere in $Q$.
To show that $\varphi\geq0$ almost everywhere in $Q$, we exploit the fact that $F(r)=+\infty$ for all $r<0$, together with \eqref{eq:regularized_potential}-\eqref{eq:regularized_potential2}.
We know that
\begin{align}
    \int_QF_\lambda(\varphi_{\lambda})\leq M_{\tau}.
    \label{eq:426_0}
\end{align}
Since $F_\lambda$ is bounded below and there is a closed interval $I\subset\mathbb{R}$ such that $\{r\in\mathbb{R}\mid F_\lambda(r)<0\}\subset I$, we arrive at
\begin{align}
    \int_{\{(x,t)\in Q\mid F_\lambda(\varphi_{\lambda}(x,t))<0\}}\bigl(-F_\lambda(\varphi_{\lambda})\bigr)\leq M\biggl|\biggl\{(x,t)\in Q\mid F_\lambda(\varphi_{\lambda}(x,t))<0\biggr\}\biggr|\leq M.
    \label{eq:426_1}
\end{align}
Observe now that
\begin{align*}
    \int_Q|F_\lambda(\varphi_{\lambda})|=\int_{\{(x,t)\in Q\mid F_\lambda(\varphi_{\lambda}(x,t))\geq0\}}F_\lambda(\varphi_{\lambda})+\int_{\{(x,t)\in Q\mid F_\lambda(\varphi_{\lambda}(x,t))<0\}}\bigl(-F_\lambda(\varphi_{\lambda})\bigr).
\end{align*}
Then, thanks to \eqref{eq:426_0} and \eqref{eq:426_1}, we conclude that
\begin{align*}
    \int_Q|F_\lambda(\varphi_{\lambda})|\leq M_{\tau}.
\end{align*}
In particular, setting
\begin{align*}
    E:=\biggl\{(x,t)\in Q\mid\varphi(x,t)<0\biggr\}\subset Q,
\end{align*}
we have
\begin{align}
    \int_E|F_{\lambda}(\varphi_{\lambda})|\leq\int_Q|F_\lambda(\varphi_{\lambda})|\leq M_{\tau}.
    \label{eq:contradd_2be_reached}
\end{align}
By contradiction, assume that the set $E$ has positive measure.
Then, in view of the strong convergence of $\{\varphi_{\lambda}\}_\lambda$, for almost every $(x,t)\in E$ we can say that also $\varphi_{\lambda}(x,t)<0$, at least for $\lambda$ small enough (up to a subsequence).
Recalling that $F_\lambda(r)=-\frac{r^3}{\lambda}-\frac{r^3}{3}$ for all $r<0$, we have
\begin{align*}
    |F_\lambda(\varphi_{\lambda}(x,t))|=\frac{|\varphi_{\lambda}(x,t)|^3}{3\lambda}(\lambda+3),
\end{align*}
with the right-hand side diverging to infinity as $\lambda$ approaches zero.
Thus, by the Fatou Lemma we get
\begin{align*}
    \liminf_{\lambda\to0^+}\int_E |F_\lambda(\varphi_{\lambda})|=+\infty,
\end{align*}
which contradicts \eqref{eq:contradd_2be_reached}, hence $|E|=0$.
This proves the lower bound $\varphi\geq0$ almost everywhere in $Q$, thus \eqref{eq:MP_varphi} holds.

Finally, recalling \ref{eq:assumption3} and assuming \eqref{eq:sigma0_in01}-\eqref{eq:h2_nonnegative}, arguing as in \cite[Subsection $3.4$]{Scarpa2021}, we get \eqref{eq:MP_sigma}.
This ends the proof of Theorem \ref{thm:existence_weak_tau>0}.

\section{Further results in the viscous case} \label{sec:moreanalysis_with_tau>0}

In these three subsections, we prove the regularity result, the validity of the strict separation property, and the continuous dependence estimate.

\subsection{Proof of Theorem \ref{thm:strong_solutions_tau>0}.} \label{subsec:additional_regularity}

In this subsection, we prove Theorem \ref{thm:strong_solutions_tau>0}, proceeding with a formal approach. A rigorous justification can be constructed using the approximation scheme described in Subsection \ref{subsec:The finite dimensional regularized problem}.
First, we address the regularity of the nutrient, aiming to prove \eqref{eq:regularity_sigma_bis}.
In this direction, we multiply \eqref{eq:sigma} by $-\Delta\sigma$ and integrate over \modLuca{$Q_t$} for any $t\in[0,T]$.
Applying the Young inequality and recalling \ref{eq:assumption2}-\ref{eq:assumption3}, we get the estimate
\begin{equation*}
    \frac{1}{2}\|\nabla\sigma(t)\|^2+\frac{1}{2}\int_{Q_t}|\Delta\sigma|^2+B\int_{Q_t}|\nabla\sigma|^2\leq M\biggl(1+\frac{1}{2}\|\sigma_0\|_V^2+\|\sigma\|_{L^2(0,T;H)}^2\biggr)\leq M.
\end{equation*}
Using elliptic regularity theory, this implies that
\begin{equation}
    \sigma\in L^\infty(0,T;V)\cap L^2(0,T;W).
    \label{eq:regularity_sigma_step}
\end{equation}
\\
Further, we reali\sz e that, by comparison in \eqref{eq:sigma} (see \ref{eq:assumption2}-\ref{eq:assumption3}), thanks to \eqref{eq:regularity_sigma_step}, we get $\sigma\in H^1(0,T;H)$.

Let us now prove \eqref{eq:additional_reg_varphi}-\eqref{eq:additional_reg_mu}.
We multiply \eqref{eq:varphi} by $\partial_t\mu$ and the time derivative of \eqref{eq:mu} by $-\partial_t\varphi$. Then we integrate over $\Omega$ and add the resulting
identities. This gives, for almost every $t\in(0,T)$, 
\begin{align*}
   \frac 12  \diff{}{t}\biggl(\|\nabla&\mu\|^2+{\tau} \|\partial_t\varphi\|^2\biggr)+\int_\Omega (F''(\varphi)+a)|\partial_t\varphi|^2\\
    =&\int_\Omega S(\varphi,\sigma)\partial_t\mu+\int_\Omega(J*\partial_t\varphi)\partial_t\varphi+\chi\int_{Q_t}\partial_t\sigma\partial_t\varphi.
\end{align*}
Setting
\begin{align}
    \mathcal{J}:=\frac{1}{2}\|\nabla\mu\|^2-\int_\Omega S(\varphi,\sigma)\mu+\frac{\tau}{2}\|\partial_t\varphi\|^2,
    \label{eq:J_functional}
\end{align}
and using \ref{ass:kernel} and \eqref{eq:property_of_F_second_0}, arguing as in \cite{Rocca2023}, we  find that
\begin{align}
    \diff{\mathcal{J}}{t}
    \leq-\int_\Omega \partial_t S(\varphi,\sigma)\mu+\int_\Omega(J*\partial_t\varphi)\partial_t\varphi+\chi\int_\Omega\partial_t\sigma\partial_t\varphi.
    \label{eq:estimate_for_J}
\end{align}
On the other hand, using  \ref{eq:assumption1}, \eqref{eq:MP_varphi}, the Poincaré--Wirtinger inequality, equation \eqref{eq:mu} and the Young inequality, we get
\begin{align}
\begin{split}
    -\int_\Omega S(\varphi,\sigma)\mu&=-\int_\Omega S(\varphi,\sigma)(\mu-\mu_\Omega)-\int_\Omega S(\varphi,\sigma)\mu_\Omega\\
    &\geq-M(\|\mu-\mu_\Omega\|+|\mu_\Omega|)\\
    &\geq-M(\|\nabla\mu\|+\|\partial_t\varphi\|_{L^1(\Omega)}+\|\varphi\|_{L^1(\Omega)}+\|\sigma\|_{L^1(\Omega)}+\|F'(\varphi)\|_{L^1(\Omega)})\\
    &\geq-\frac{1}{8}\|\nabla\mu\|^2-\delta_1\|\partial_t\varphi\|^2-\ov M(1+\|\varphi\|_{L^1(\Omega)}+\|\sigma\|_{L^1(\Omega)}+\|F'(\varphi)\|_{L^1(\Omega)}),
    \label{eq:source_in_J}
\end{split}
\end{align}
for some constant $\ov M>0$ and for all $\delta_1>0$.
Recalling now \eqref{eq:MP_varphi} and \eqref{eq:separation_for_mean}, from \eqref{eq:mz_applied_lambda} we infer
\begin{align}
    \|F'(\varphi)\|_{L^1(\Omega)}\leq C_1\int_\Omega F'(\varphi)(\varphi-\varphi_\Omega)+C_2|\Omega|,
    \label{eq:MZ_applied}
\end{align}
almost everywhere in $(0,t)$, for any $t\in (0,T]$.
Thanks to \eqref{eq:mu} and the fact that $\int_{\Omega}\mu_\Omega(\varphi-\varphi_\Omega)=0$, we can rewrite the right-hand side of \eqref{eq:MZ_applied} to obtain
\begin{align}
\begin{split}
    & \|F'(\varphi)\|_{L^1(\Omega)}\leq C_1\int_\Omega (\mu-\mu_\Omega)(\varphi-\varphi_\Omega)\\
    &\quad 
    +M\biggl(1-\tau\int_\Omega\partial_t\varphi(\varphi-\varphi_\Omega) 
    -\int_\Omega(a\varphi-J*\varphi)(\varphi-\varphi_\Omega)+\int_\Omega\sigma(\varphi-\varphi_\Omega)\biggr).
    \label{eq:54_4}
\end{split}
\end{align}
Recalling \eqref{eq:application_of_PW}, we have
\begin{align*}
    C_1\int_\Omega(\mu-\mu_\Omega)(\varphi-\varphi_\Omega)\leq \frac{1}{8 \ov M}\int_\Omega|\nabla\mu|^2+M\int_\Omega|\varphi|^2.
\end{align*}
The other terms on the right-hand side of \eqref{eq:54_4} can be estimated by the Young inequality, whence
\begin{align}
    \ov M\|F'(\varphi)\|_{L^1(\Omega)}\leq\frac{1}{8}\|\nabla\mu\|^2+\delta_2\|\partial_t\varphi\|^2+M\bigl(1+\|\varphi\|^2+\|\sigma\|^2\bigr),
    \label{eq:estimate_F_prime}
\end{align}
for all $\delta_2>0$.
Thus, taking into account \eqref{eq:J_functional}, \eqref{eq:source_in_J} and \eqref{eq:estimate_F_prime}, we observe that  
\begin{align*}
    \mathcal{J}\geq\frac{1}{4}\|\nabla\mu\|^2+\biggl(\frac{\tau}{2}-\delta_1-\delta_2\biggr)\|\partial_t\varphi\|^2-M.
\end{align*}
It is then possible to select $\delta_1,\delta_2>0$ in such a way that
\begin{align}
    M_\tau(\mathcal{J}+1)\geq\|\nabla\mu\|^2+\|\partial_t\varphi\|^2,
    \label{eq:J_bigger_than}
\end{align}
for a constant $M_\tau>0$ which depends on $\tau>0$.
Let us go back to \eqref{eq:estimate_for_J}.
The first term on the right-hand side can be bounded exploiting \eqref{eq:grad_h1_bdd} and the chain rule in Sobolev spaces, while the last two terms on the right-hand can be readily bounded by means of the Hölder inequality, the Young convolution inequality and the Young inequality.
We arrive at
\begin{align*}
    \diff{\mathcal{J}}{t}\leq M\bigl(1+\|\partial_t\varphi\|^2+\|\partial_t\sigma\|^2+\|\mu\|^2\bigr),
\end{align*}
whence
\begin{align*}
    \diff{\mathcal{J}}{t}\leq M_\tau\bigl(1+\mathcal{J}+\|\varphi\|^2+\|\partial_t\sigma\|^2+\|\mu\|^2\bigr),
\end{align*}
in view of \eqref{eq:J_bigger_than}.
Integrating in time and exploiting some previous bounds, we get
\begin{align*}
    \mathcal{J}(t)\leq\mathcal{J}(0)+M_\tau\int_0^t (1+\mathcal{J}(\tau))d\tau.
\end{align*}
We now claim that, from the assumptions on the initial data, it follows that $\mathcal{J}(0)\leq M$.
Then, the Gronwall Lemma allows to conclude that $\mathcal{J}(t)\leq M_\tau$ for almost every $t\in(0,T)$, i.e., by \eqref{eq:J_bigger_than}, that \eqref{eq:additional_reg_varphi} and \eqref{eq:additional_reg_mu} hold.
We have
\begin{align*}
    \mathcal{J}(0)\leq M\bigl(\|\varphi_0\|^2+\|\mu_0\|_V^2+\|\sigma_0\|^2+\|\varphi_{0}'\|^2\bigr),
\end{align*}
where $\varphi_{0}'$ and $\mu_0$ are defined by
\begin{align}
   & \varphi_{0}'= \Delta\mu_0  + S(\varphi_0,\sigma_0),\quad
    \mu_{0}=\tau\varphi_{0}'+a\varphi_{0}-J\ast\varphi_{0}+F'(\varphi_{0})-\chi\sigma_{0}, \label{eq:initial_velocity}
\end{align}
almost everywhere in $\Omega$.
Indeed, since $\varphi_0$ and $\sigma_0$ are known data, we can rewrite \eqref{eq:initial_velocity} as an elliptic reaction-diffusion problem for $\varphi_0'$, namely
\begin{align*}
    \varphi_{0}'-\tau\Delta\varphi_{0}'=\widetilde{S}(\varphi_0,\sigma_0)\quad\text{a.e. in }\Omega,
\end{align*}
where $\widetilde{S}(\varphi_0,\sigma_0):=S(\varphi_0,\sigma_0)+\Delta(a\varphi_{0}-J\ast\varphi_{0}+F'(\varphi_{0})-\chi\sigma_{0})$ is a known space-dependent function.
Since
\begin{align*}
    \|\mu_0\|_V\leq\tau\|\varphi'_0\|_V+\|a\varphi_{0}-J\ast\varphi_{0}+F'(\varphi_{0})-\chi\sigma_{0}\|_V,
\end{align*}
by a standard elliptic argument, we infer that
\begin{align*}
    \mathcal{J}(0)\leq M\bigl(1+\|\varphi_0\|_V^2+\|F'(\varphi_0)\|_V^2+\|\sigma_0\|_V^2\bigr).
\end{align*}
Then, the claim follows in view of \eqref{eq:assumptions_ics} and \eqref{eq:assumptions_ics_enhanced}.

Finally, integrating equation \eqref{eq:mu}  
with respect to time and taking the gradient of both sides of the resulting equation, we get
\begin{equation*}  
\tau\nabla \varphi = \tau\nabla\varphi_0 
-\int_0^t (\nabla(a\varphi)-\nabla F'(\varphi) - \nabla \mu -\chi\nabla\sigma - \nabla J*\varphi)\quad\text{a.e. in } \Omega,
\end{equation*}
for all $t\in (0,T]$. Thus, we
also have $\varphi \in \mathcal{C}^0([0,T];V)$ and the proof of Theorem \ref{thm:strong_solutions_tau>0} is finished.

\subsection{Proof of Theorem \ref{thm:separation_property_varphi_tau>0}} \label{subsec:separation_property}
To prove the separation property \eqref{eq:separation_property_varphi}, we draw upon an argument presented in \cite{Colli2021} (see Subsection $3.1$ therein) and adapt the proof to the single-well potential case.
We start by rewriting equation \eqref{eq:varphi} as follows
\begin{align*}
        -\Delta\mu=S(\varphi,\sigma)-\partial_t\varphi,
        \quad\text{a.e. in }Q.
\end{align*}
The right-hand side is bounded in $L^\infty(0,T;H)$ in view of \ref{eq:assumption1}, \eqref{eq:MP_varphi} and \eqref{eq:additional_reg_varphi}.
Then, a standard argument implies that $\mu\in L^\infty(0,T;W)\subset L^\infty(Q)$.
Let us now rewrite \eqref{eq:mu} as an ODE with respect to the phase variable $\varphi$:
\begin{align}    \tau\partial_t\varphi+a\varphi+F'(\varphi)=\mu+\chi\sigma+J*\varphi\quad\text{a.e. in } Q,
    \label{eq:mu_sep_property}
\end{align}
where the right-hand side is bounded in $L^\infty(Q)$. 
Since $F\in\mathcal{C}^\infty([0,1))$ with $\lim_{r\to1^-}F'(r)=+\infty$, we get the existence of a constant $\delta\in(0,\bar{\delta}]$ such that
\begin{align}
    F'(r)-\|\mu+\chi\sigma+J*\varphi\|_{L^\infty(Q)}\geq0, \quad\forall r\in(1-\delta,1).
    \label{eq:potential_sep_property}
\end{align}
Multiplying \eqref{eq:mu_sep_property} by $(\varphi-(1-\delta))_+$ and integrating over \modLuca{$Q_t$} for any $t\in[0,T]$, owing to \eqref{eq:potential_sep_property}, we obtain
\begin{align*}
    \frac{\tau}{2}\|(\varphi(t)-(1-\delta))_+\|^2+a_*\int_{Q_t}\varphi(\varphi-(1-\delta))_+\leq\frac{\tau}{2}\|(\varphi_0-(1-\delta))_+\|^2.
\end{align*}
Recalling that $a_*\int_{Q_t}\varphi(\varphi-(1-\delta))_+\geq0$, and (see \eqref{eq:separation_for_ic}) $
    \frac{\tau}{2}\|(\varphi_0-(1-\delta))_+\|^2=0, 
$
we end up with $(\varphi-(1-\delta))_+=0$ almost everywhere in $Q$.

This ends the proof of Theorem \ref{thm:separation_property_varphi_tau>0}.

\subsection{Proof of Theorem \ref{thm:continuous_dependence_tau>0}} \label{subsec:stability_estimate}

Let us set
\begin{equation*}
    \varphi:=\varphi_1-\varphi_2,\quad\mu:=\mu_1-\mu_2,\quad\sigma:=\sigma_1-\sigma_2, \quad \varphi_0:=\varphi_0^1-\varphi_0^2,\quad\sigma_0:=\sigma_0^1-\sigma_0^2.
\end{equation*}
It is \modLuca{straightforward} to check that $(\varphi,\mu,\sigma)$ solves the following problem
\begin{align}
    &\partial_t\varphi-\Delta\mu=S(\varphi_1,\sigma_1)-S(\varphi_2,\sigma_2)&\text{ in }Q,
    \label{eq:uniqueness_phi} \\
    &\mu=\tau\partial_t\varphi+a\varphi-J\ast\varphi+F'(\varphi_1)-F'(\varphi_2)-\chi\sigma&\text{ in }Q,
    \label{eq:uniqueness_mu} \\
    &\partial_t\sigma-\Delta\sigma+B\sigma+C\sigma {\hhd}(\varphi_1)+C\sigma_2({\hhd}(\varphi_1)-{\hhd}(\varphi_2))=0&\text{ in }Q,
    \label{eq:uniqueness_sigma} \\
    &\partial_{\boldsymbol{n}}\mu=\partial_{\boldsymbol{n}}\sigma=0&\text{ on }\Sigma, \label{eq:uniqueness_bcs} \\
    &\varphi(0)=\varphi_0,\quad\sigma(0)=\sigma_0,\qquad&\text{in }\Omega. \label{eq:uniqueness_ics}
\end{align}
Let us multiply \eqref{eq:uniqueness_phi} by $\mathcal{N}(\varphi-\varphi_\Omega)$, \eqref{eq:uniqueness_mu} by $-(\varphi-\varphi_\Omega)$ and \eqref{eq:uniqueness_sigma} by $\sigma$.
Adding up the resulting equalities and integrating over $Q_t$ for every $t\in[0,T]$, we obtain
\begin{align*}
    \int_0^t\langle\partial_t&\varphi,\mathcal{N}(\varphi-\varphi_\Omega)\rangle+\int_0^t\langle-\Delta\mu,\mathcal{N}(\varphi-\varphi_\Omega)\rangle-\int_{Q_t}\mu(\varphi-\varphi_\Omega)+\frac{1}{2}\|\sigma(t)\|^2\\
    &+\int_{Q_t}\bigl( a|\varphi|^2+(F'(\varphi_1)-F'(\varphi_2))\varphi\bigr)+\frac{\tau}{2}\|\varphi(t)\|^2+\int_{Q_t}|\nabla\sigma|^2+B\int_{Q_t}|\sigma|^2\\
    =&C\int_{Q_t}|\sigma|^2\hhd(\varphi_1)+C\int_{Q_t}\sigma\sigma_2({\hhd}(\varphi_2)-{\hhd}(\varphi_1))+\frac{1}{2}\|\sigma_0\|^2+\frac{\tau}{2}\|\varphi_0\|^2\\
    &+\int_{Q_t}(S(\varphi_1,\sigma_1)-S(\varphi_2,\sigma_2))\mathcal{N}(\varphi-\varphi_\Omega)+\int_{Q_t}\bigl(J*\varphi+\chi\sigma\bigr)\varphi\\
    &+\int_{Q_t}\bigl(\tau\partial_t\varphi+a\varphi-J\ast\varphi+F'(\varphi_1)-F'(\varphi_2)-\chi\sigma\bigr)\varphi_\Omega.
\end{align*}
Recalling that $\mathcal{N}(\varphi-\varphi_\Omega)$ has zero mean value and that $-\Delta\mu=-\Delta(\mu-\mu_\Omega)$, exploiting \eqref{eq:property_of_N_1} and \ref{eq:assumption2}, from the above identity, we infer the following estimate:

\begin{align}
\begin{split}
    \frac{1}{2}\|(\varphi&-\varphi_\Omega)(t)\|^2_{V'}+\frac{1}{2}\|\sigma(t)\|^2+\int_{Q_t}|\nabla\sigma|^2+B\int_{Q_t}|\sigma|^2\\
    &+\int_{Q_t}\bigl( a|\varphi|^2+(F'(\varphi_1)-F'(\varphi_2))\varphi\bigr)+\frac{\tau}{2}\|\varphi(t)\|^2\\
    \leq&C\|{\hhd}\|_{L^\infty(\mathbb{R})}\int_{Q_t}|\sigma|^2+C\int_{Q_t}\sigma\sigma_2({\hhd}(\varphi_2)-{\hhd}(\varphi_1))+\int_{Q_t}\bigl(J*\varphi+\chi\sigma\bigr)\varphi\\
    &+\frac{1}{2}\|\varphi_0-(\varphi_0)_\Omega\|^2_{V'}+\frac{\tau}{2}\|\varphi_0\|^2+\frac{1}{2}\|\sigma_0\|^2+\int_{Q_t}(S(\varphi_1,\sigma_1)-S(\varphi_2,\sigma_2))\mathcal{N}(\varphi-\varphi_\Omega)\\
    &+\int_{Q_t}\bigl(\tau\partial_t\varphi+a\varphi-J\ast\varphi+F'(\varphi_1)-F'(\varphi_2)-\chi\sigma\bigr)\varphi_\Omega.
    \label{eq:uniqueness_step1} 
\end{split}
\end{align}
Observe that
\begin{align*}
\begin{split}
    \int_{Q_t}\bigl( a|\varphi|^2+(F'(\varphi_1)-F'(\varphi_2))\varphi\bigr)&\geq a_*\int_{Q_t}|\varphi|^2+\int_{Q_t}(F'(\varphi_1)-F'(\varphi_2))\varphi\\
    &=\int_{Q_t}\biggl(a_*+\frac{F'(\varphi_1)-F'(\varphi_2)}{\varphi_1-\varphi_2}\biggr)|\varphi|^2\geq c_0\int_{Q_t}|\varphi|^2,
\end{split}
\end{align*}
where the last inequality follows from \eqref{eq:property_of_F_second_0} (recall that $F\in\mathcal{C}^2([0,1))$).
Then, taking into account \eqref{eq:MP_sigma} and applying the Young inequality, from \eqref{eq:uniqueness_step1} we get
\begin{align}
\begin{split}
    \frac{1}{2}\|(\varphi&-\varphi_\Omega)(t)\|^2_{V'}+\frac{1}{2}\|\sigma(t)\|^2+\int_{Q_t}|\nabla\sigma|^2+B\int_{Q_t}|\sigma|^2+c_0\int_{Q_t}|\varphi|^2+\frac{\tau}{2}\|\varphi(t)\|^2\\
    \leq&\frac{1}{2}\|\varphi_0-(\varphi_0)_\Omega\|^2_{V'}+\frac{\tau}{2}\|\varphi_0\|^2+\frac{1}{2}\|\sigma_0\|^2\\
    &+M\biggl(\int_{Q_t}|\varphi|^2+\int_{Q_t}|\sigma|^2\biggr)+\int_{Q_t}(S(\varphi_1,\sigma_1)-S(\varphi_2,\sigma_2))\mathcal{N}(\varphi-\varphi_\Omega)\\
    &+\int_{Q_t}\bigl(\tau\partial_t\varphi+a\varphi-J\ast\varphi+F'(\varphi_1)-F'(\varphi_2)-\chi\sigma\bigr)\varphi_\Omega.
    \label{eq:uniqueness_step2}
\end{split}
\end{align}
On the other hand, we have
\begin{align*}
    \int_{Q_t}(a\varphi-J\ast\varphi)\varphi_\Omega=\int_0^t\varphi_\Omega\int_\Omega(a\varphi-J\ast\varphi)=0.
\end{align*}
Moreover, since $\varphi_1$ and $\varphi_2$ satisfy \eqref{eq:separation_property_varphi}, in view of the Lipschitz continuity of $F$ on the interval $[0,1-\delta]$ for all $\delta\in(0,1)$, thanks to the Jensen inequality, we get
\begin{align*}
    \int_{Q_t}(F'(\varphi_1)-F'(\varphi_2))\varphi_\Omega\leq M_\tau\int_{Q_t}|\varphi|^2.
\end{align*}
The chemotactic term can be estimated by means of the Young and Jensen inequalities. This gives
\begin{align*}
    -\chi\int_{Q_t}\sigma\varphi_\Omega\leq M\biggl(\int_{Q_t}|\varphi|^2+\int_{Q_t}|\sigma|^2\biggr).
\end{align*}
Finally, we observe that
\begin{align*}
    \tau\int_{Q_t}\partial_t\varphi\varphi_\Omega=\tau\int_0^t\varphi_\Omega(s)\int_\Omega\partial_t\varphi(s)\dd s=\tau|\Omega|\int_0^t\varphi_\Omega(s)\varphi'_\Omega(s)\dd s.
\end{align*}
To control this last term, we follow \cite{Grasselli2024} and multiply \eqref{eq:uniqueness_phi} by $|\Omega|^{-1}$ to obtain
\begin{align*}
    \varphi'_\Omega=-m\varphi_\Omega+\frac{1}{|\Omega|}\int_\Omega({\hhu}(\varphi_1,\sigma_1)-{\hhu}(\varphi_2,\sigma_2)),
\end{align*}
for almost every $t\in(0,T)$.
Multiplying the above equation by $\varphi_\Omega$ and then appealing to the Young inequality, \eqref{eq:grad_h1_bdd} and the Jensen inequality, one gets
\begin{align*}
    \varphi_\Omega \, \varphi'_\Omega\leq M(\|\varphi\|^2+\|\sigma\|^2).
\end{align*}
Taking into account all the previous remarks, \eqref{eq:uniqueness_step2} becomes
\begin{align*}
\begin{split}
    \frac{1}{2}\|(\varphi&-\varphi_\Omega)(t)\|^2_{V'}+\frac{1}{2}\|\sigma(t)\|^2+\int_{Q_t}|\nabla\sigma|^2+c_0\int_{Q_t}|\varphi|^2+\frac{\tau}{2}\|\varphi(t)\|^2\\
    \leq&\frac{1}{2}\|\varphi_0-(\varphi_0)_\Omega\|^2_{V'}+\frac{\tau}{2}\|\varphi_0\|^2+\frac{1}{2}\|\sigma_0\|^2\\
    &+M_\tau\biggl(\int_{Q_t}|\varphi|^2+\int_{Q_t}|\sigma|^2\biggr)+\int_{Q_t}(S(\varphi_1,\sigma_1)-S(\varphi_2,\sigma_2))\mathcal{N}(\varphi-\varphi_\Omega).
\end{split}
\end{align*}
We are left to estimate the last term on the right-hand side.
Notice that $\mathcal{N}(\varphi-\varphi_\Omega)$ is mean-free, so we can subtract to the quantity $S(\varphi_1,\sigma_1)-S(\varphi_2,\sigma_2)$ its mean.
Thus, appealing to \eqref{eq:property_of_N_2}, the Young inequality, the Jensen inequality and \eqref{eq:grad_h1_bdd}, we arrive at
\begin{align*}
    & \frac{1}{2}\|(\varphi-\varphi_\Omega)(t)\|^2_{V'}+\frac{1}{2}\|\sigma(t)\|^2+\int_{Q_t}|\nabla\sigma|^2+c_0\int_{Q_t}|\varphi|^2+\frac{\tau}{2}\|\varphi(t)\|^2+m\|\varphi-\varphi_\Omega\|^2_{L^2(0,T;V')}\\
    & \quad \leq M_\tau\biggl(\int_{Q_t}|\varphi|^2+\int_{Q_t}|\sigma|^2\biggr)+\frac{1}{2}\|\varphi_0-(\varphi_0)_\Omega\|^2_{V'}+\frac{\tau}{2}\|\varphi_0\|^2+\frac{1}{2}\|\sigma_0\|^2
    \\ & \qquad 
    +\int_0^t\|(\varphi-\varphi_\Omega)(s)\|^2_{V'}\dd s,
\end{align*}
and an application of the Gronwall Lemma gives the proof.

\section{Proof of Theorem \ref{thm:asymptotic_tau_to_0}} \label{sec:asymptotics_as_tau_to_0}
We now proceed formally, as the argument can be rigorously justified using an approximation scheme similar to that employed in Subsection \ref{subsec:The finite dimensional regularized problem}.
Before proceeding with the uniform estimates in $\tau$, observe that a weak solution $(\varphi_\tau,\mu_\tau,\sigma_\tau)$ to problem \eqref{eq:varphi}-\eqref{eq:ics} in the sense of Definition \ref{def:weak_sol} exists, with $\varphi_\tau$ satisfying \eqref{eq:MP_varphi} and \eqref{eq:separation_for_mean}, where the constant $\delta$ in \eqref{eq:separation_for_mean} is independent of $\tau$.

We now set out to retrace the energy estimate taking care to estimate some terms in such a way as to obtain $\tau$-independent estimates, and then send it to the limit $\tau \to 0$.
Thus, we consider \eqref{eq:53_0}.
To this, we add 
\eqref{eq:53_1}, without performing the estimate \eqref{eq:grad_grad_estimate} to control the term involving $\chi$.
This leads to (cf. also \eqref{eq:property_of_F_second_0})
\begin{align*}
	\non
	    & \int_{Q_t}|\nabla\mu_{\tau}|^2+\tau\int_{Q_t}|\partial_t\varphi_{\tau}|^2+\epsilon\|\varphi_{\tau}(t)\|^2+\frac{1}{2}\|\sigma_{\tau}(t)\|^2+\int_{Q_t}|\nabla\sigma_{\tau}|^2
	    \\
& \non \qquad 
 +\frac{1}{2}\|\varphi_\tau(t)\|^2+\frac{\tau}{2}\|\nabla\varphi_\tau(t)\|^2+\frac{c_0}{2}\int_{Q_t}|\nabla\varphi_\tau|^2	    
	    \\
   &
   \non \quad  \leq  a^*\|\varphi_{0,\tau}\|^2+\frac{1}{2}\|\sigma_{0,\tau}\|^2
+\|\varphi_{0,\tau}\|^2+\frac{\tau}{2}\|\varphi_{0,\tau}\|_V^2   
   +\int_\Omega F(\varphi_{0,\tau})
    \\ & \qquad \non 
    +M\biggl(1+\int_{Q_t}|\varphi_\tau|^2 +\int_{Q_t}|\sigma_\tau|^2 \biggr)
    \\
    &\qquad 
    +\chi\int_{Q_t}\nabla\sigma_\tau\cdot\nabla\varphi_\tau
    +\chi\int_{Q_t}\sigma_\tau\partial_t\varphi_\tau+\int_{Q_t}S(\varphi_\tau,\sigma_\tau)\mu_\tau.
\end{align*}
We are left to control the last three terms on the \rhs. To bound the last one, we can argue as in Subsection \ref{subsec:Uniform estimates} (see \eqref{eq:fondamentale2}).
Then, to estimate the two terms that are multiplied by $\chi$, we first integrate by parts, further appealing to equation \eqref{eq:sigma} and to the Young inequality.

Namely, we observe that for all $\delta_1,\delta_2>0$ it holds
\begin{align*}
\begin{split}
    \chi\int_{Q_t}\sigma_{\tau}&\partial_t\varphi_{\tau}+\chi\int_{Q_t}\nabla\varphi_{\tau}\cdot\nabla\sigma_{\tau}\\
    =&-\chi\int_0^t\langle\partial_t\sigma_{\tau}(s),\varphi_{\tau}(s)\rangle\,ds+\chi\int_\Omega\sigma_{\tau}(t)\varphi_{\tau}(t)\\
    &-\chi\int_\Omega\sigma_{0,\tau}\varphi_{0,\tau}+\chi\int_{Q_t}\nabla\varphi_{\tau}\cdot\nabla\sigma_{\tau} \\
    \leq& M\biggl(1+\|\varphi_{0,\tau}\|^2+\|\sigma_{0,\tau}\|^2+\int_{Q_t}|\varphi_{\tau}|^2+\int_{Q_t}|\sigma_{\tau}|^2\biggr)\\
    &+\delta_1\chi^2\|\varphi_{\tau}(t)\|^2+\frac{1}{4\delta_1}\|\sigma_{\tau}(t)\|^2+4\delta_2\chi^2\int_{Q_t}|\nabla\varphi_{\tau}|^2+\frac{1}{4\delta_2}\int_{Q_t}|\nabla\sigma_{\tau}|^2.
\end{split}
\end{align*}
Thus, we get
\begin{align*}
    \int_{Q_t}&|\nabla\mu_\tau|^2+\frac{\tau}{2}\int_{Q_t}|\partial_t\varphi_\tau|^2+\epsilon\|\varphi_\tau(t)\|^2+\frac{1}{2}\|\sigma_\tau(t)\|^2+\int_{Q_t}|\nabla\sigma_\tau|^2\\
    &+\frac{1}{2}\|\varphi_\tau(t)\|^2+\frac{\tau}{2}\|\nabla\varphi_\tau(t)\|^2+\frac{c_0}{2}\int_{Q_t}|\nabla\varphi_\tau|^2\\
    \leq& M\bigl(1+\|\varphi_{0,\tau}\|^2+\tau\|\varphi_{0,\tau}\|_V^2+\|\sigma_{0,\tau}\|^2+\|F(\varphi_{0,\tau})\|_{L^1(\Omega)}\an{}\bigl)\\
    &+M\biggl(\int_{Q_t}|\varphi_\tau|^2+\int_{Q_t}|\sigma_\tau|^2\biggr)+\delta_1\chi^2\|\varphi_{\tau}(t)\|^2\\
    &+\frac{1}{4\delta_1}\|\sigma_{\tau}(t)\|^2+4\delta_2\chi^2\int_{Q_t}|\nabla\varphi_{\tau}|^2+\frac{1}{4\delta_2}\int_{Q_t}|\nabla\sigma_{\tau}|^2.
\end{align*}
The terms referring to the initial conditions are all bounded due to assumptions \eqref{eq:lastthm_ics1} and \eqref{eq:lastthm_ics2}, and the last four terms can be absorbed into the left-hand side provided we choose $\delta_1,\delta_2>0$ such that
\begin{align*}
\begin{split}
\begin{cases}
    \epsilon+\frac{1}{2}-\delta_1\chi^2>0,\\
    \frac{1}{2}-\frac{1}{4\delta_1}>0,\\
    \frac{c_0}{2}-4\delta_2\chi^2>0,\\
    1-\frac{1}{4\delta_2}>0,
\end{cases}
\end{split}
\end{align*}
which is possible thanks to assumption \eqref{eq:tau_critic_0}.
Finally, the Gronwall Lemma yields
\begin{align*}
    \|\varphi&_\tau\|_{L^\infty(0,T;H)\cap L^2(0,T;V)}+\|\mu_\tau\|_{L^2(0,T;V)}\\
    &+\|\sigma_\tau\|_{L^\infty(0,T;H)\cap L^2(0,T;V)}+\tau^\frac{1}{2}\|\varphi_\tau\|_{H^1(0,T;H)\cap L^\infty(0,T;V)}\leq M,
\end{align*}
for a constant $M$ that is independent of $\tau$.
We can now conclude the proof by using a standard argument based on the estimate obtained above, which is independent of $ \tau $. More precisely, repeating a procedure similar to the one in Subsection \ref{subsec:Passage to the limit}, it is possible to let $ \tau $ go to zero along a suitable subsequence and finding a limit triple which solves problem \eqref{eq:varphi}-\eqref{eq:ics} with $\tau=0$ in the sense of Definition \ref{def:weak_sol}.

\section*{Acknowledgements}
MG and AS gratefully acknowledge partial support
from the MIUR-PRIN Grant 2020F3NCPX ``Mathematics for industry 4.0 (Math4I4)'', support by MUR, grant Dipartimento di Eccellenza 2023-2027
and their affiliation to the GNAMPA (Gruppo Nazionale per l'Analisi Matematica, la Probabilit\`a e le loro Applicazioni) of INdAM (Isti\-tuto Nazionale di Alta Matematica). AS also acknowledges support from the Alexander von Humboldt Foundation, and LM was partially supported by EPSRC (Engineering and Physical Sciences Research Council).


\end{document}